\documentclass[11pt,reqno,a4paper]{amsart}

\usepackage{anysize}
\marginsize{3.5cm}{3.5cm}{2.5cm}{2.5cm}
\usepackage[english]{babel}
\usepackage{amsmath}
\usepackage{amsfonts,amssymb}
\usepackage{enumerate}
\usepackage{mathrsfs}
\usepackage{amsthm}
\usepackage{tikz}
\usetikzlibrary{arrows}
\usepackage{mathrsfs}
\usepackage{caption}

\usepackage{hyperref} 
\hypersetup{
    colorlinks=true,       
    linkcolor=red,          
    citecolor=blue,        
    filecolor=magenta,      
    urlcolor=cyan           
}

\theoremstyle{plain}
\newtheorem{theo}{Theorem}[section]
\newtheorem*{theo*}{Theorem}
\newtheorem{prop}[theo]{Proposition}
\newtheorem{lemm}[theo]{Lemma}

\newtheorem{defi}[theo]{Definition}

\theoremstyle{definition}
\newtheorem{rema}[theo]{Remark}
\newtheorem*{rema*}{Remark}

\newtheorem*{nota*}{Notation}

\DeclareMathOperator{\cn}{div}

\DeclareMathOperator{\diff}{d}

\DeclareSymbolFont{pletters}{OT1}{cmr}{m}{sl}
\DeclareMathSymbol{s}{\mathalpha}{pletters}{`s}
\renewcommand{\a}{\alpha}

\renewcommand{\b}{\beta}
\def\aeta{\arrowvert_{y=\eta}}
\def\ai{\sqrt{1+\eta_x^2}}

\def\ba{\begin{align}}

\def\bad{\begin{aligned}}
\def\be{\begin{equation}}
\def\blA{\bigl\lVert}

\def\brA{\bigr\rVert}

\def\cutoff{\varphi}
\def\curv{\px\left(\frac{\eta_x}{\ai}\right)}

\def\defn{\mathrel{:=}}

\def\dsigma{\diff \! \sigma}
\def\dt{\diff \! t}
\def\dx{\diff \! x}

\def\dxdt{\dx \dt}
\def\dy{\diff \! y}

\def\dydx{\diff \! y \diff \! x}
\def\dydxdt{\diff \! y \diff \! x \diff \! t}

\def\e{\eqref}
\def\ea{\end{align}}
\def\ead{\end{aligned}}

\def\ee{\end{equation}}
\def\eps{\varepsilon}

\def\fracE{\frac{\diff\! \mathcal{E}}{\dt}}

\def\fract{\frac{\diff}{\dt}}

\def\Hr{\mathcal{H}}

\def\intL{\int_{-L}^L}
\def\intT{\int_0^T}
\def\intTL{\int_0^T\int_{-L}^L}

\def\la{\left\vert}
\def\lA{\left\Vert}
\def\le{\leq}

\def\mez{\frac{1}{2}}

\def\px{\partial_x}

\def\ra{\right\vert}

\def\rA{\right\Vert}

\def\ST{\Big\arrowvert_0^T}

\def\tdm{\frac{3}{2}}

\def\uq{\frac{1}{4}}

\def\xR{\mathbb{R}}

\def\xT{\mathbb{T}}
\def\xZ{\mathbb{Z}}

\numberwithin{equation}{section}

\title{Stabilization of the water-wave equations with surface tension}
\author[Thomas Alazard]{Thomas Alazard \\ CNRS \& \'Ecole Normale Sup\'erieure Paris-Saclay}
\def\notina[#1]#2{\begingroup\def\thefootnote{\fnsymbol{footnote}}\footnote[#1]{#2}\endgroup}

\begin{document}

\begin{abstract}
This paper is devoted to the stabilization of the water-wave equations with surface tension through of an external pressure acting on a small part of the free surface. It is proved that the energy decays to zero exponentially in time, provided that the external pressure is given by the normal component of the velocity 
at the free surface multiplied by an appropriate cut-off function.
\end{abstract}

\maketitle

\section{Introduction}\label{S:2}

Consider the incompressible Euler equations for a potential flow\notina[0]{This work is 
partly supported by the grant ``ANA\'E'' ANR-13-BS01-0010-03.}
in a fluid domain located between with a free surface, two vertical walls and a flat bottom, which is at time $t$ of the form
$$
\Omega(t)= \{\, (x,y)\in [0,L] \times \xR\, : \, -h<y<\eta(t,x)\,\},
$$
where $L$ is the length of the basin, $h$ is its depth and $\eta$, the free surface elevation, is an unknown function. 
We assume that the flow is irrotational so that 
the eulerian velocity field is the gradient of a potential function $\phi=\phi(t,x,y)$ satisfying
\begin{equation}\label{t5}
\begin{aligned}
&\Delta_{x,y}\phi=0\quad\text{in }\Omega(t),\\ 
&\partial_{t} \phi +\mez \la \nabla_{x,y}\phi\ra^2 +P +g y = 0 \quad\text{in }\Omega(t),\\
&\phi_y =0 \quad \text{on }y=-h,\\
&\phi_x=0\quad\text{on }x=0 \text{ or }x=L,
\end{aligned}
\end{equation}
where 
$P\colon \Omega\rightarrow\xR$ is the pressure, 
$g>0$ is the acceleration of gravity, 
$\nabla_{x,y}=(\partial_x,\partial_y)$ and $\Delta_{x,y}=\partial_x^2+\partial_y^2$. 
Partial differentiations in space will be denoted by suffixes so that $\phi_x=\partial_x\phi$ and $\phi_y=\partial_y \phi$.

The water-wave equations are then given by two boundary conditions on the free surface: 
the classical kinematic boundary condition,
describing the deformations of the domain,
\be\label{t8}
\partial_{t} \eta = \sqrt{1+\eta_x^2}\, \phi_n \arrowvert_{y=\eta}=\phi_y(t,x,\eta(t,x))-
\eta_x(t,x)\phi_x(t,x,\eta(t,x)),
\ee
together with an equation expressing the balance of forces across the free surface:
\be\label{t9}
P\arrowvert_{y=\eta}=P_{ext}-\kappa H(\eta),
\ee
where $P_{ext}=P_{ext}(t,x)$ is the evaluation of the external pressure at the free surface, $\kappa\ge 0$ is the coefficient of surface tension 
and $H(\eta)$ is the curvature of the free surface:
\begin{equation*}
H(\eta)=\px\left( \frac{\eta_x}{\sqrt{1+\eta_x^2}}\right).
\end{equation*}
We are concerned with the problem with surface tension and assume that $\kappa>0$.

Introduce the energy $\mathcal{E}=\mathcal{E}(t)$, defined by
\be\label{t2ini}
\mathcal{E}=\frac{g}{2}\int_{0}^L\eta^2\dx+\kappa \int_{0}^L \left(\sqrt{1+\eta_x^2}-1\right)\dx
+\mez\int_{0}^L\int_{-h}^{\eta(t,x)}\la \nabla_{x,y}\phi \ra^2\dydx.
\ee
This is the sum of the gravitational potential energy, 
a surface energy due to stretching of the surface and the kinetic energy. 
Recall that the energy is conserved when there is no external pressure, 
which means that if $P_{ext}=0$ then 
$\mathcal{E}(t)=\mathcal{E}(0)$ for all time. 
The stabilization problem for the water-wave equations consists in finding 
a pressure law, relating $P_{ext}$ to the unknown $(\eta,\psi)$, such that:
\begin{enumerate}[i)]
\item $\mathcal{E}$ is decreasing and converges to zero exponentially in time;
\item $\partial_x P_{ext}$ is supported inside a small subset of $[0,L]$. 
\end{enumerate}

The stabilization problem 
corresponds to an important issue 
in the numerical analysis of water waves, namely the problem of 
damping outgoing waves in an absorbing zone surrounding the computational boundary. 
There is a huge literature about the absorption of the water-wave energy 
in a numerical wave tank and we refer the 
reader to \cite{Bonnefoy2005,CBS1993,Clamond2005,Clement1996,Clement1999,Duclos2001,Ducrozet2007,Grilli1997,JPCKRA,JKR} 
and references therein. A popular choice
is to assume that $P_{ext}=\chi \partial_t \eta$ for some cut-off function $\chi\ge 0$ 
supported\footnote{We consider a damping located near $x=L$ only, and not also near $x=0$, since we imagine that water waves are generated near $x=0$. 
To do so, one could use also the variations of an external pressure. 
The latter problem is related to control theory and some references are given below.} in $ [L-\delta,L]$ 
for some $\delta>0$. To explain this choice, we start by recalling that, 
as observed by Zakharov~\cite{Zakharov1968}, the equations can be written in hamiltonian form. To do so, Zakharov 
introduced  
$\psi(t,x)=\phi(t,x,\eta(t,x))$, observed that 
the energy can be written as a function of $(\eta,\psi)$ and verified that
$$
\frac{\partial\eta}{\partial t}=\frac{\delta \mathcal{E}}{\delta \psi},\qquad
\frac{\partial\psi}{\partial t}=-\frac{\delta \mathcal{E}}{\delta \eta}-P_{ext}.
$$
Using these equations, one deduces that
\be\label{p3}
\fracE=\int \left(\frac{\delta \mathcal{E}}{\delta \eta}\frac{\partial\eta}{\partial t}+\frac{\delta \mathcal{E}}{\delta \psi}\frac{\partial\psi}{\partial t}\right)\dx
=-\int \frac{\partial\eta}{\partial t} P_{ext}\dx,
\ee
and hence, if $P_{ext}=\chi \partial_t \eta$ with $\chi\ge 0$, we deduce that $\diff \!\mathcal{E}/\dt\le 0$. 
It is thus easily seen that the energy decays. However, it is much more complicated to prove that the energy converges exponentially to zero. 
To study this problem, we first need to pause to clarify the question since, in general, 
solutions of the water-wave equations do not exist globally in time (they might 
blow-up in finite time, see \cite{CCFGGS,CS2}). 
Our goal is in fact to prove that there exists a constant $C$ 
such that, if a regular solution exists on the time interval $[0,T]$, then
\be\label{p1}
\mathcal{E}(T)\le \frac{C}{T} \mathcal{E}(0).
\ee
Since the equation is invariant by translation in time, one can iterate this inequality. Consequently, if the solution exists on time intervals of size $nT_0$ with 
$T_0\ge 2C$, then $\mathcal{E}(nT_0)\le 2^{-n}\mathcal{E}(0)$, which is the desired exponential decay. 

We studied a similar problem in \cite{A-Stab-WW} 
for the case without surface tension. Assuming that $\kappa=0$ and defining $P_{ext}$ by
\be\label{p2}
\partial_x P_{ext}=\chi(x)\int_{-h}^{\eta(t,x)}\phi_x(t,x,y)\dy,
\ee
where $\chi\ge 0$ is a cut-off function, we proved in \cite{A-Stab-WW} an inequality of the form
\be\label{p5}
\mathcal{E}(T)\le \frac{C(N)}{\sqrt{T}}\mathcal{E}(0),
\ee
where the constant $C(N)$ depends on the frequency 
localization\footnote{The quantity $N$ measuring the frequency localization is of the ratio of two Sobolev norms. 
One can think of the ratio $\lA u\rA_{H^1}/\lA u\rA_{L^2}$, which is proportional to $N$ for a typical function 
oscillating at frequency $N$, like $u(x)=\cos(2\pi Nx/L)$.} of the solution $(\eta,\psi)$. 
The fact that this constant must depend on the frequency localization can be easily understood by considering 
the linearized equations (see \cite{A-Stab-WW}). Somewhat surprisingly, it was possible to prove this result 
under some mild smallness condition on the solution which allows to consider truly nonlinear regimes. However, 
the main drawback is that one cannot control the frequency localization of the solutions in terms of the initial data, 
unless one makes strong smallness assumptions on these initial data (see \S$5$ in \cite{Boundary}). In this paper, we will prove that one can exploit surface tension to improve the stabilization of water waves in three directions:
\begin{enumerate}[i)]
\item the main improvement is that, under mild smallness assumptions on $\eta$, 
one can prove that $\mathcal{E}(T)\le (C/T)\mathcal{E}(0)$ for some constant $C$ depending only on the physical parameters 
$g,\kappa,h,L$. In particular, $C$ is independent of the frequency localization. 
Another key point is that we will give an {\em explicit} bound for $C$. 
As a consequence, under explicit mild smallness assumptions on $\eta$, 
it is possible to use the induction argument mentioned above to obtain exponential decay 
of the energy. By combining this result with the local 
controllability result for small data proved by Alazard, Baldi and Han-Kwan in \cite{ABHK}, 
this in turn implies a controllability result for non small data in large time (following an argument due to 
Dehman, Lebeau and Zuazua~\cite{DLZ2003}).
\item The second important improvement is that we are able to consider the case when $P_{ext}=\chi\partial_t\eta$ instead of assuming that $P_{ext}$ is given by \e{p2}. This is important for applications since, as mentioned above, this feedback law is widely used in the numerical analysis of water waves. 
\item Eventually we obtain an inequality of the form \e{p1} with a factor $1/T$ instead of the factor $1/\sqrt{T}$ of \e{p5}. 
\end{enumerate}

\section{Main result}

Following Zakharov~\cite{Zakharov1968} and Craig--Sulem~\cite{CrSu} (see also \cite{Bertinoro}), 
one reduces the water-wave equations to a system on the free surface. 
To do so, remarking that the velocity potential $\phi$ 
is harmonic and satisfies a Neumann boundary condition on the walls and the bottom, 
we see that $\phi$ is fully determined by its evaluation at the free surface. 
We thus work below with
$$
\psi(t,x)\defn\phi(t,x,\eta(t,x)),
$$
and introduce the Dirichlet to Neumann operator, denoted by~$G(\eta)$, relating 
$\psi$ to the normal derivative of the potential by 
$$
G(\eta)\psi=\sqrt{1+\eta_x^2}\, \phi_n \arrowvert_{y=\eta}=(\phi_y-\eta_x \phi_x)\aeta.
$$
Then, it follows from \e{t8} that $\partial_t \eta=G(\eta)\psi$ and \e{t5}--\e{t9} implies that
\begin{align}
&\partial_t \psi+g\eta +N(\eta)\psi-\kappa H(\eta)=-P_{ext},\quad\text{where}\notag\\
&
N(\eta)\psi=\mathcal{N} \big\arrowvert_{y=\eta}\quad\text{with}\quad
\mathcal{N}=
\mez\phi_x^2-\mez \phi_y^2+\eta_x \phi_x \phi_y,\label{t90}
\end{align}
where recall that $H(\eta)=\px\left(\eta_x / \sqrt{1+\eta_x^2}\right)$. 
One can check that 
$$
N(\eta)\psi=\mez \psi_x^2-\mez \frac{(G(\eta)\psi+\eta_x\psi_x)^2}{1+\eta_x^2}.
$$

With the above notations, the water-wave system reads
\be\label{systemT}
\left\{
\begin{aligned}
&\partial_t \eta=G(\eta)\psi,\\
&\partial_t \psi+g\eta +N(\eta)\psi-\kappa H(\eta)=-P_{ext}.
\end{aligned}
\right.
\ee
In this paper we consider solutions $(\eta,\psi)$ of \e{systemT}Ê
which are periodic in $x$, with period $2L$ and we further assume 
that these solutions are even in $x$. 
As explained below, the reason why we are making these symmetry assumptions is that this setting 
corresponds to the setting discussed in the introduction. That is the 
case where the water waves are propagating in a bounded container with length $L$, 
having vertical walls located on $x=0$ and $x=L$, so that at time $t$ the fluid domain is 
\be\label{j1}
\Omega(t)= \{\, (x,y)\in [0,L] \times \xR\, : \, -h<y<\eta(t,x)\,\}.
\ee

\begin{defi}[Regular solutions]\label{D1}
Let $L>0$, $T>0$. Denote by $\xT$ the torus $\xR/(2L \xZ)$ and, for any $\sigma\in\xR$, 
denote by $H^\sigma(\xT)$ the Sobolev space of order $\sigma$ of $2L$-periodic functions. We say that $(\eta,\psi,P_{ext})$ is 
a regular solution of \e{systemT} defined on $[0,T]$ provided that 
the following three conditions are satisfied:
\begin{enumerate}[i)]
\item there exists a real number $s>5/2$ such that
\begin{align*}
&\eta\in C^0([0,T];H^{s+\mez}(\xT)),\quad \psi\in C^0([0,T];H^{s}(\xT)),\\
&P_{ext}\in C^0([0,T];H^{s-1}(\xT));
\end{align*}
\item $\eta$, $\psi$ and $P_{ext}$ are even in $x$;
\item for all $t$ in $[0,T]$,  
\be\label{t9b}
\inf_{x\in [-L,L]}\eta(t,x)\ge -\frac{h}{2},\quad \int_{-L}^{L}\eta(t,x) \dx=0,\quad 
\int_{-L}^L P_{ext}(t,x)\dx=0.
\ee
\end{enumerate}
\end{defi}
\begin{rema}
$i)$ We can make without loss of generality an 
assumption on the mean value of $\eta$ since it is a conserved quantity. 

$ii)$ Let us recall from \cite{ABZ4} an argument which shows 
that this setting allows us to consider water waves in a bounded container. 
Consider a regular solution $(\eta,\psi,P_{ext})$ of \e{systemT} defined on $[0,T]$. 
The assumption that $\eta$ and $\psi$ are even and periodic with period $2L$ implies that 
$\phi$ is also even and periodic in $x$. We deduce that 
$\phi_x$ is odd in $x$ and hence 
$\phi_x(t,x,y)=0$ whenever $x=0$ or $x=L$. 
In particular, the normal component of the velocity satisfies the solid wall boundary condition (that 
is $v\cdot n=0$) on both the bottom and the lateral walls ($\{x=0\}$ or $\{x=L\}$). 
Consequently, as illustrated below, we obtain a solution to the problem 
in a bounded container by considering the restrictions of $\eta$ to $0\le x\le L$ 
and the restriction of $\phi$ to $\Omega$ given by~\e{j1} (see Figure~\ref{fig2} below).
\end{rema}

\def\a{1.7671}
\def\deuxa{3.5343}
\def\b{6.4795}
\def\c{8.2467}
\def\cbis{8.2967}

\begin{center}
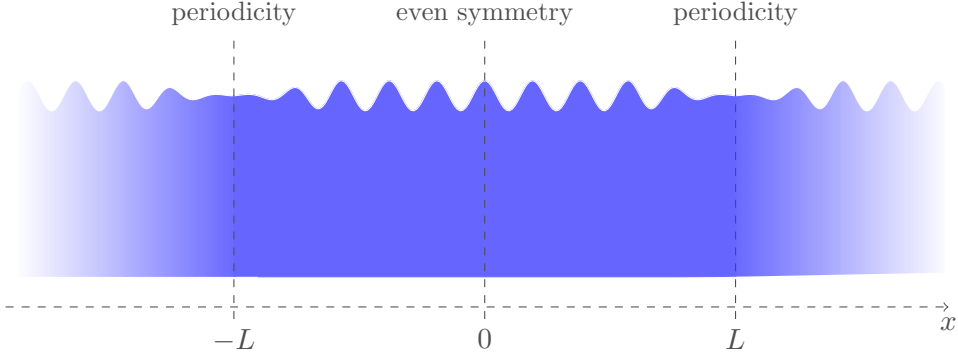
\begin{figure}[h]
\begin{tikzpicture}[scale=0.8,samples=100]

\shadedraw  [right color=blue!60!white, left color=blue!00!white, color=white, shift={(-\cbis,0)}] 
(-\deuxa,-3) -- (-\deuxa,0) --  
plot [domain=-\a:0,shift={(-\a,0)}] ({\x},{0.25*cos(8*\x r)}) -- 
 plot [domain=0:\a,shift={(-\a,0)}] ({\x},{0.25*exp(-\x*\x)*cos(8*\x r)}) -- (0.1,0) -- (0.1,-3) ;

\filldraw [color=white, shift={(-\cbis,0)}] 
(-3.7,-3) -- (-3.7,0.3) -- (-3.5,0.3) -- (-3.5,-3);

\shadedraw  [left color=blue!60!white, right color=blue!00!white, color=white,shift={(\c,0)}] 
(-\cbis,-3) -- (-\cbis,0) --  plot [domain=-\a:0,shift={(-\b,0)}] ({\x},{0.25*exp(-\x*\x)*cos(8*\x r)}) -- 
plot [domain=0:\a,shift={(-\b,0)}] ({\x},{0.25*cos(8*\x r)}) -- (-4.7124,-3);

\filldraw [color=white, shift={(\c,0)}] 
(-4.8,-3) -- (-4.8,0.3) -- (-4.7,0.3) -- (-4.7,-3);

 \filldraw  [color=blue!60!white] 
(-\c,-3) -- (-\c,0) --  plot [domain=-\a:0,shift={(-\b,0)}] ({\x},{0.25*exp(-\x*\x)*cos(8*\x r)}) -- 
plot [domain=0:4.7124,shift={(-\b,0)}] ({\x},{0.25*cos(8*\x r)}) -- 
 plot [domain=0:\a,shift={(-\a,0)}] ({\x},{0.25*exp(-\x*\x)*cos(8*\x r)}) -- (0,-3) -- (-7.8540,-3)  ;

 \draw  [color=white] 
(-\c,0) --  plot [domain=-\a:0,shift={(-\b,0)}] ({\x},{0.25*exp(-\x*\x)*cos(8*\x r)}) -- 
plot [domain=0:4.7124,shift={(-\b,0)}] ({\x},{0.25*cos(8*\x r)}) -- 
 plot [domain=0:\a,shift={(-\a,0)}] ({\x},{0.25*exp(-\x*\x)*cos(8*\x r)}) ;
 
 \filldraw  [color=white]
 (-10,-3.2) -- (10,-3.2) -- (10,-2.8) ;

\draw [color=black!70!white,dashed] (-4.1234,-3.7) node [below] {$0$}-- (-4.1234,1) node [above] {\small{even symmetry}};
\draw [color=black!70!white,dashed] (0,-3.7) node [below] {$L$}-- (0,1) node [above] {\small{periodicity}};
\draw [color=black!70!white,dashed] (-\c,-3.7) node [below] {$-L$}-- (-\c,1) node [above] {\small{periodicity}};
\draw [color=black!70!white,->,dashed] (-12.0,-3.5) -- (3.5,-3.5) node [below] {$x$};
\end{tikzpicture}
\caption{Periodic and even solutions.}\label{fig2}
\end{figure}
\end{center}

We denote by $\Hr$ the energy
\be\label{t2}
\mathcal{H}=\frac{g}{2}\int_{-L}^L\eta^2\dx+\kappa \int_{-L}^L \left(\sqrt{1+\eta_x^2}-1\right)\dx
+\mez\int_{-L}^L\int_{-h}^{\eta(t,x)}\la \nabla_{x,y}\phi \ra^2\dydx.
\ee
\begin{rema}\label{R:24}
Since $\eta$ and $\phi$ are even and periodic in $x$ with period $2L$, we deduce that 
$\Hr=2\mathcal{E}$ where $\mathcal{E}$ is as defined 
in the introduction (see \e{t2ini}), that is
$$
\mathcal{E}=\frac{g}{2}\int_{0}^L\eta^2\dx+\kappa \int_{0}^L \left(\sqrt{1+\eta_x^2}-1\right)\dx
+\mez\int_{0}^L\int_{-h}^{\eta(t,x)}\la \nabla_{x,y}\phi \ra^2\dydx.
$$
\end{rema}

We now define the pressure law. We find convenient to impose that the mean value\footnote{In the introduction, we mentioned that the pressure law used in the literature is 
simply $P_{ext}=\chi\partial_t\eta$. Let us give two arguments showing that one can add a time dependent function to the pressure. The first observation is that,  
since ${\fract \int\eta\dx=0}$ (conservation of the volume), for any function $F=F(t)$ depending only on time we have 
$\int F \partial_t \eta \dx=0$ so ${\int \frac{\partial\eta}{\partial t} P_{ext}\dx=\int \frac{\partial\eta}{\partial t}( P_{ext}+F)\dx}$. 
Consequently, the argument given above in the introduction 
to deduce that the energy decays (see \e{p3}) still applies when $P_{ext}$ is replaced by 
$P_{ext}+F(t)$. The second observation is that, when $P_{ext}$ is given by \e{p4}, the pressure is not 
supported inside a small subset of $[-L,L]$ but this is not an obstruction since we only require that its spatial 
derivative is supported inside a small subset of $[-L,L]$. Indeed, 
the only quantity which is physically meaningful is the velocity which is the gradient 
of the velocity potential. In other words, we could also work with $P_{ext}=\chi\partial_t \eta$ by 
modifying the Bernouilli's constant in the equation for $\psi$.} of $P_{ext}$ 
is zero, so we set
\be\label{p4}
P_{ext}(t,x)=\chi(x)\partial_t \eta (t,x)-\frac{1}{2L}\int_{-L}^L \chi(x)\partial_t \eta (t,x)\dx,
\ee
for some cut-off function $\chi\ge 0$ to be determined. In order to exploit a cancellation in the proof, we choose 
$\chi$ as follows. 

\begin{defi}\label{defi:P}
Fix $\delta>0$ and consider a $2L$-periodic $C^\infty$ cut-off function $\cutoff$, satisfying $0\le \cutoff\le 1$ 
and such that 
$$
\cutoff(x)=\cutoff(-x),\quad x\cutoff'(x)\le 0 \text{ for }x\in [-L,L],\quad  
\cutoff(x)=\left\{\begin{aligned}&1 \text{ if }x\in [0,L-\delta],\\ &0 \text{ if }x\in \left[L-\frac{\delta}{2},L\right].
\end{aligned}\right.
$$
We successively set  
$$
m(x)=x\cutoff(x),
$$
and
$$
\chi(x)=1-m_x(x).
$$
\end{defi}

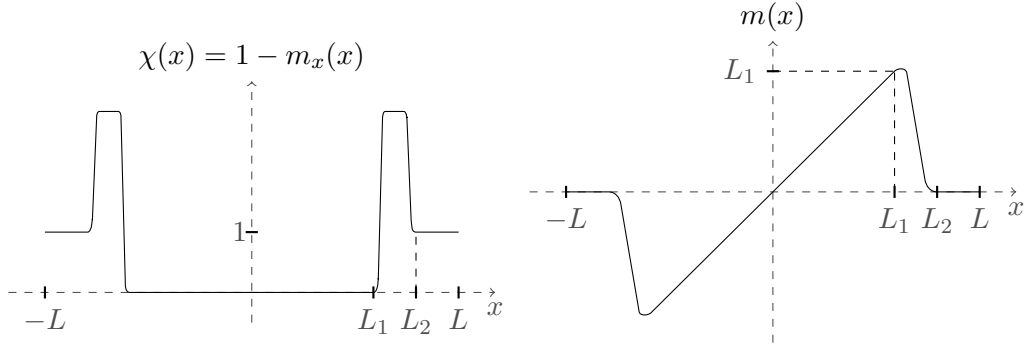
\begin{figure}[h!]
\begin{minipage}[b]{0.40\linewidth}
\begin{tikzpicture}[scale=0.8,samples=100]
\draw [color=black!70!white,->,dashed] (-4,0) -- (4,-0) node [below] {$x$};
\draw [color=black!70!white,->,dashed] (0,-0.5) -- (0,3.5) node [color=black,above] {$\chi(x)=1-m_x(x)$};
\draw[thick] (3.4,0.1) -- (3.4,-0.1) node[color=black!70!white,below] {$L$};
\draw[thick] (-3.4,0.1) -- (-3.4,-0.1) node[color=black!70!white,below] {$-L$};
\draw[thick] (-0.1,1) -- (0.1,1) node[color=black!70!white,left] {$1$};
\draw[thick] (2,0.1) -- (2,-0.1) node[color=black!70!white,below] {$L_1$};
\draw [dashed] (2.7,0) -- (2.7,1);
\draw[thick] (2.7,0.1) -- (2.7,-0.1) node[color=black!70!white,below] {$L_2$};
\draw (-3.4,1) -- (-2.7,1) to [out=0,in=90] (-2.62,1.2) -- (-2.55,2.9) to [out=90,in=180] (-2.5,3) -- (-2.2,3) to [out=0,in=90] (-2.15,2.9) -- (-2.08,0.2) 
to [out=-90,in=180] (-2,0) -- 
(0,0) -- (2,0) to [out=0,in=90] (2.08,0.219) -- (2.15,2.9) to [out=80,in=180] (2.2,3) -- (2.5,3) to [out=0,in=-80] (2.55,2.9) -- (2.62,1.2) 
to  [out=-80,in=180] (2.7,1) -- (3.4,1) ;
\end{tikzpicture}
\hspace{10mm}
\end{minipage}
\hspace{10mm}
\begin{minipage}[b]{0.40\linewidth}
\begin{tikzpicture}[scale=0.8,samples=100]
\draw [color=black!70!white,->,dashed] (-4,0) -- (4,-0) node [below] {$x$};
\draw [color=black!70!white,->,dashed] (0,-2.5) -- (0,2.5) node [color=black,above] {$m(x)$};
\draw[thick] (3.4,0.1) -- (3.4,-0.1) node[color=black!70!white,below] {$L$};
\draw[thick] (-3.4,0.1) -- (-3.4,-0.1) node[color=black!70!white,below] {$-L$};
\draw (-3.4,0) -- (-2.7,0) to [out=0,in=100] (-2.5,-0.2) -- (-2.2,-1.98) to [out=-85,in=-135] (-2,-2) -- 
(0,0) -- (2,2) to [out=45,in=105] (2.2,1.98) -- (2.5,0.2)  to [out=-80,in=180] (2.7,0) -- (3.4,0);
\draw [dashed] (2,0) -- (2,2);
\draw [dashed] (0,2) -- (2,2);
\draw[thick] (0.1,2) -- (-0.1,2) node[color=black!70!white,left] {$L_1$};
\draw[thick] (2,0.1) -- (2,-0.1) node[color=black!70!white,below] {$L_1$};
\draw[thick] (2.7,0.1) -- (2.7,-0.1) node[color=black!70!white,below] {$L_2$};
\end{tikzpicture}
\hspace{10mm}
\end{minipage}
\caption{The functions $\chi$ and $m$, here $L_1=L-\delta$ and $L_2=L-\delta/2$.\label{fig1}}
\end{figure}
\begin{rema}
Let $m$ be as the given by Figure~\ref{fig1} and define $\cutoff$ by $\cutoff(x)=m(x)/x$. 
One immediately see that
$$
0\le \varphi\le 1,\quad \cutoff(x)=\cutoff(-x),\quad  
\cutoff(x)=\left\{\begin{aligned}&1 \text{ if }x\in [0,L-\delta],\\ &0 \text{ if }x\in \left[L-\frac{\delta}{2},L\right].
\end{aligned}\right.
$$
It remains to check that $x\cutoff'(x)\le 0 \text{ for }x\in [-L,L]$. 
We have 
$x\cutoff'(x)=f(x)/x^2$ with $f(x)=xm'(x)-m(x)$. Set $L'=(L_1+L_2)/2$, so that $m$ is concave on $[0,L']$ and convex on $[L',L]$. 
Observe that $f(0)=0$ and $f'(x)=xm''(x)\le 0$ for $x\in [0,L']$ since $m$ is concave there, 
which proves that $f(x)\le 0$ on $[0,L']$. On the other hand, for $x\ge L'$, one has $m'(x)\le 0$ and $m(x)\ge 0$ so one obviously get $f(x)\le 0$. This 
proves that $x\cutoff'(x)\le 0$ on $[0,L]$. Since $x\cutoff'(x)$ is even in $x$, we deduce that $x\cutoff'(x)\le 0$ on $[-L,L]$.
\end{rema}

Here is our main result. Recall that $\kappa$ is the surface tension coefficient, $g$ is the acceleration 
of gravity, $h$ is the depth of the fluid domain at rest and $L$ the width of the basin (half of the period).

\begin{theo}\label{T:2}Let $m,\chi$ be as given by Definition~$\ref{defi:P}$. 
Assume that $\kappa>0$ and
\be\label{t1000}
\kappa \sup_{[-L,L]}m_{xx}(x)^2\le g.
\ee
Then 
there exists a positive constant $C$, depending only on the physical parameters 
$\kappa,g,h,L$, such that the following result holds. Let $T>0$ and consider a regular solution $(\eta,\psi,P_{ext})$ 
defined on the time interval $[0,T]$ with $P_{ext}$ as given by~$\e{p4}$. Set 
$$
\rho(t,x)=(m(x)-x) \eta_x(t,x) +\frac{9}{4}\eta(t,x) -\mez m_x(x)\eta(t,x).
$$
If, for all $(t,x)\in [0,T]\times [-L,L]$,
\begin{alignat*}{2}
&i)\qquad &&\rho(t,x) \ge -\frac{h}{4},\qquad  \la \rho_x(t,x)\ra< \frac 1 4,\\
&ii) \qquad &&\int_{-L}^L (1-m_x(x))\eta(t,x) \dx \le \frac{hL}{3},\qquad 
\la m_x(x)\ra\la \eta_x(t,x)\ra^2\le 2,\quad \la \eta_x(t,x)\ra\le 1,
\end{alignat*}
then one has the estimate
$$
\Hr(T)\le \frac{C}{T}\Hr(0).
$$
\end{theo}
\begin{rema}
$i)$ One can choose $m$ as in Figure~\ref{fig1} so that $\la m_{xx}\ra \le 8\delta^{-1}+32L\delta^{-2}$. Then \e{t1000} is 
satisfied provided that $\kappa\le (8\delta^{-1}+32L\delta^{-2})^{-2}g$. 

$ii)$ An important remark is that the constant $C$ can be given by an explicit formula in terms of $\kappa,g,h,L$. 
Namely, $C$ is given by $C=8K$ where $K$ is as defined in \e{d20} (taking $\lambda=1$ in \e{d20}).
\end{rema}

To conclude this introduction, we outline the main steps of the proof and also 
comment on related questions about the study of the Cauchy problem or the controllability of the water-wave equations.

\subsection{About the Cauchy problem}
As already mentioned, Theorem~\ref{T:2} 
implies that the energy converges exponentially to zero in the following sense: let $T_0=2C$ where 
$C$ is as given by Theorem~\ref{T:2} and consider a regular solution $(\eta,\psi,P_{ext})$ defined on the time interval $[0,nT_0]$, then one has
$$
\Hr(nT_0)\le 2^{-n}\Hr(0).
$$
This remark raises a question about the Cauchy problem, namely to prove that 
one can obtain regular solutions on arbitrary large time intervals. 
We will not study this problem in this paper because 
it involves different tools. Here, we only  
state a result, which will be proved in a separate paper, implying that regular solutions 
exist on large time intervals, under a smallness assumption on the initial data. 

Consider the equations 
\be\label{systemT2}
\left\{
\begin{aligned}
&\partial_t \eta=G(\eta)\psi,\\
&\partial_t \psi+g\eta +N(\eta)\psi-\kappa H(\eta)=-P_{ext},
\end{aligned}
\right.
\ee
where $P_{ext}$ is given by 
\be\label{systemT3}
P_{ext}(t,x)=\chi(x)\partial_t \eta -\frac{1}{2L}\int_{-L}^L \chi(x)\partial_t \eta (t,x)\dx,
\ee
for some function $\chi\ge 0$ (the following result holds for any function $\chi$, 
and in particular we do not need here to assume that $\chi$ can be written as $1-m_x$ for some multiplier $m$ of the form $m(x)=x\cutoff(x)$). 
Consider now an initial data $(\eta_0,\psi_0)$ in $H^{s+\mez}(\xT)\times H^s(\xT)$ for some $s>5/2$, 
$2L$-periodic and even, and set
$$
m=\lA \eta_0\rA_{H^{s+\mez}(\xT)}+\lA \psi_0\rA_{H^{s}(\xT)}.
$$
Then the Cauchy problem for the system \e{systemT2}--\e{systemT3} has a unique regular solution (which means 
that $(\eta,\psi,P_{ext})$ satisfies the three properties of Definition~$\ref{D1}$) 
on the time interval $[0,c_*/m]$ for some positive constant $c_*$ depending only on $s$. 
Moreover, one has the estimate
$$
\sup_{t\in [0,c_*/m]}\Big\{\lA \eta(t)\rA_{H^{s+\mez}(\xT)}+\lA \psi(t)\rA_{H^{s}(\xT)}\Big\}\le 2m.
$$

\subsection{Application to the controllability of the water-wave equations}
There are many results about the controllability or the stabilization of linear or nonlinear equations 
describing water waves in some asymptotic regimes like Benjamin-Ono, KdV or nonlinear Schr\"odinger equations 
(we refer the reader to the book of Coron~\cite{Coron}). 
However, one cannot easily adapt these studies to the water-wave system~\e{systemT2} since 
the latter system is quasi-linear (instead of semi-linear) and since  
it is a pseudo-differential system, 
involving the Dirichlet-Neumann operator which is 
nonlocal and also depends nonlinearly on the unknown. 
The first results about the possible applications of control theory to the water-wave equations 
are due to Reid and Russell \cite{ReidRussell1985} and Reid~\cite{Reid1986,Reid1995}, who studied 
the linearized equations at the origin (see also Miller~\cite{Miller2012}, 
Lissy~\cite{Lissy2014} and Biccari~\cite{Biccari} 
for other control results 
about dispersive equations involving a fractional Laplacian). 
Alazard, Baldi and Han-Kwan proved in \cite{ABHK} the first controllability result for 
the nonlinear water-wave equations with surface tension, namely a controllability result in arbitrarily small time, 
under a smallness assumption on the size of the data. 

\begin{theo}[Alazard, Baldi, Han-Kwan, from \cite{ABHK}]\label{TABHK}
Assume that $\kappa>0$. Let $T>0$ and consider a non-empty open subset $\omega\subset \xT$. 
There exist $\sigma$ large enough and a positive constant 
$\eps_c$ small enough such that, for any two pairs of functions $(\eta_{in},\psi_{in})$,  
$(\eta_{final},\psi_{final})$ in $H^{\sigma+\mez}(\xT)\times H^\sigma(\xT)$ satisfying 
$$
\lA \eta_{in}\rA_{H^{\sigma+\mez}}+\lA \psi_{in}\rA_{H^\sigma}<\eps_c,\quad 
\lA \eta_{final}\rA_{H^{\sigma+\mez}}+\lA \psi_{final}\rA_{H^\sigma}<\eps_c,
$$
and such that $\int \eta_{in}\dx =\int \eta_{final}\dx=0$, 
there exists $P_{ext}$ in $C^0([0,T];H^{\sigma}(\xT))$, supported in 
$[0,T]\times \omega$, such that the Cauchy problem for \e{systemT} has a unique solution
$$
(\eta,\psi)\in C^0([0,T];H^{\sigma+\mez}(\xT)\times H^{\sigma}(\xT)), 
$$
and the solution $(\eta,\psi)$ satisfies 
$(\eta\arrowvert_{t=T},\psi\arrowvert_{t=T})=(\eta_{final},\psi_{final})$.
\end{theo}

The smallness condition in \cite{ABHK} comes 
from various technical arguments and it is possibly a strong assumption. 
In sharp contrast, Theorem~\ref{T:2} holds under explicit and mild smallness conditions on $\eta$. 
One can combine both results to obtain a controllability result for larger initial data in large time. 
This strategy was used by Dehman-Lebeau-Zuazua~\cite{DLZ2003} and Laurent \cite{Laurent_2010a,Laurent_2010b} 
for other wave equations. 
The idea is to proceed in two steps: firstly one steers 
the solution close to zero by taking as a control
the stabilization term so that it is possible, in a second step, to use the local controllability near zero to steer the solution to zero. 
This proves a null controllability result. 
By using another standard argument, exploiting the fact that the equation is reversible in time, 
one deduces the wanted controllability result from this null controllability result. 

\usetikzlibrary{decorations}
\usepgflibrary{decorations.pathmorphing}

\begin{center}
\begin{figure}[h]
\begin{tikzpicture}[scale=0.8,samples=100]
\node at (0,3) {\textbullet};
\node at (0,3) [right]{$u_{in}$};
\node at (0,3) [left]{large data};
\node at (7,0) {\textbullet};
\node at (7,0) [below]{$0$};
\draw[color=blue,thick] plot [domain=0:5] ({\x},{3*exp(-0.3*\x)});
\draw [dashed] (5,0) -- (5,4) node[left]{{\small{stabilization}}};
\draw [dashed] (7,0) -- (7,3.5) ;
\node at (7,4) {{\small{null controllability}}};
\node at (6,3) {{\small{forward}}};
\node at (8,3) {{\small{backward}}};
\draw [dashed] (9,0) -- (9,4) node[right]{{\small{stabilization}}};
\draw[color=green,thick] decorate [decoration={random steps,segment length=1pt,amplitude=0.6pt}] 
    {(4.995,0.668) -- (7,0)};
\draw[color=green,thick] decorate [decoration={random steps,segment length=1pt,amplitude=0.6pt}] 
    {(7,0) -- (9.005,0.668)};
\draw [dashed] (0,0.668) -- (12.5,0.668) ;
\node at (0,0.668) [left]{small data};
\draw[color=blue,thick] plot [domain=9:12] ({\x},{0.57+2.3*exp((\x-12))}) node[color=black,right]{$u_{final}$};
\node at (12,2.9) {\textbullet};
\draw [thick,->] (0,0) -- (0,4) node[above]{size};
\draw[thick,->] (0,0) -- (12.5,0) node[right]{time};
\end{tikzpicture}
\caption{The Dehman-Lebeau-Zuazua strategy (cf~\cite{DLZ2003}).}
\end{figure}
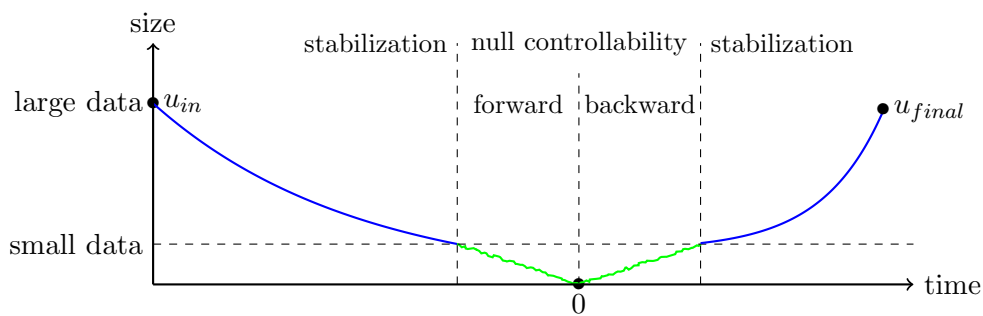
\end{center}

\subsection{Strategy of the proof}\label{S:23}
To prove Theorem~\ref{T:2}, our goal is to 
prove that there exists a positive constant $C$ such that 
\be\label{I11}
\int_0^T\Hr(t)\dt\le C\Hr(0).
\ee
Since $\Hr$ is a decreasing function, this will imply the wanted inequality
$$
\Hr(T)\le \frac{1}{T} \int_0^T\Hr(t)\dt\le \frac{C}{T}\Hr(0).
$$
The proof of \e{I11} is in two steps.

\smallbreak

\textbf{First step.} 
Here we derive an exact {\em identity} which involves the integral 
in time of the energy. This approach is not new.  It was already
performed in our previous works \cite{Boundary,A-Stab-WW} and was based on several tools: 
the multiplier technique (with the multiplier $m(x)\partial_x$ 
for some function $m$ to be determined), the Craig-Sulem-Zakharov 
reduction to an hamiltonian system on the boundary, a Pohozaev identity for the 
Dirichlet to Neumann operator, previous results about the Cauchy problem and computations 
guided by the analysis done by Benjamin and Olver of the conservation laws for water waves (cf \cite{BO}). 
We need here to adapt this analysis to the present context with surface tension, 
which requires new ideas. Indeed, compared to the case without 
surface tension, we will have to handle many remainder terms. 
In the end, instead of stating an identity (as we did for the case 
without surface tension), 
we will obtain an inequality. 

Consider a smooth 
function $m=m(x)$ with $m(L)=0$, and set
$$
\zeta=\partial_x(m\eta)+\tdm (1-m_x)\eta-\uq \eta,
$$
our first main task will be to derive the following inequality
\be\label{I12}
\begin{aligned}
&\frac 1 4 \int_0^T \Hr(t)\dt\le O+W+B-I\qquad \text{where}\\ 
I&\defn \frac{h}{4}\int_0^T\int_{-L}^L\phi_x(t,x,-h)^2\dxdt,\\
O&\defn \int_0^T\int_{-L}^L \left(\tdm(1-m_x)\widetilde{\psi}+(x-m)\psi_x\right)G(\eta)\psi\dxdt, \\
W&\defn -\int_0^T\int_{-L}^L P_{ext}\, \zeta\, \dxdt,\\
B&\defn \int_{-L}^L \zeta(0,x) \widetilde{\psi}(0,x) \dx-\int_{-L}^L \zeta(T,x) \widetilde{\psi}(T,x) \dx,
\end{aligned}
\ee
where $\widetilde{\psi}$ is equal to $\psi$ minus its average on $[-L,L]$. 
These four quantities play different roles. Their key properties are the following:
\begin{itemize}
\item $I\ge 0$ and hence \e{I12} gives a bound for the horizontal component of the velocity 
at the bottom (this plays a key role to control the velocity in terms of the pressure, see \e{last1}).
\item $W$ is the only term which involves the pressure.
\item $O$ corresponds to an {\em observation}, this means that this term depends 
only on the behavior of the solutions near $\{x=-L\}$ or $\{x=L\}$ when $m$ is as given by Definition~\ref{defi:P}. Indeed, 
$x-m$ and $1-m_x$ vanish when $x\in [-L+\delta,L-\delta]$ by definition of $m$. 
\item $B$ is not an integral in time, by contrast with the other terms and, in addition, it is easily 
estimated by $K\Hr(0)$.
\end{itemize}

\textbf{Second step.} The goal of the second step is to deduce the wanted result \e{I11} from 
the identity~\e{I12} proved in the first step. 
To do so, it is sufficient to prove that 
there exists a constant $K$ depending only on $g,\kappa,h,L$ such that 
\be\label{C13ini}
O+W+B -I \le K \Hr(0)+a \int_0^T\Hr(t)\dt \quad\text{for some}\quad a <\uq.
\ee
Indeed, by combining \e{I12} and \e{C13ini} one obtains that
$$
\int_0^T\Hr(t)\dt\le \frac{K}{1/4-a}\Hr(0),
$$
which is the wanted result~\e{I11}.
To prove \e{C13ini}, one estimates separately the terms $B,W,O$. 

The estimate of the term $B$ is easy so we merely explain how to estimate $W,O$.  
To explain the first ingredient of the proof, recall that 
the hamiltonian structure of the equation implies that
$$
\fract \Hr(t)=-\int_{-L}^L P_{ext} \partial_t \eta\dx.
$$
Since $\partial_t \eta=G(\eta)\psi$, one deduces that
$$
\intTL \chi (\partial_t\eta)^2\dxdt=\intTL \chi (G(\eta)\psi)^2\dxdt =\intTL P_{ext} \partial_t \eta \dxdt \le \Hr(0).
$$
This gives an estimate for the $L^2$-norm of $\chi\partial_t\eta$, which is the main 
contribution to the definition of $P_{ext}$ (see \e{systemT3}). 
A more tricky inequality, relying on the special choice $\chi(x)=1-m_x(x)$, is that
\be\label{chichoice}
-\intT \left( \int_{-L}^L \chi\partial_t \eta \dx\right) \intL \zeta\dxdt\le \frac{L}{g} \sup_{[-L,L]}(1-m_x)^2 \Hr(0).
\ee
By combining the above inequalities, we will be able to estimate the term $W$. 
The main difficulty is to bound the observation term~$O$, and in 
particular to estimate the contribution of
$$
\iint (x-m)\psi_x G(\eta)\psi\dxdt.
$$
Using the Cauchy-Schwarz inequality, the key point is to 
estimate the $L^2$-norm of $\psi_x$ in terms of the two quantities which are under control, that is 
the integral of $\chi (G(\eta)\psi)^2$ and the positive term $I$. 
In this direction, we will prove the following inequality, which is of independent interest: 
there exists a constant $A$, depending only on $\lA \eta_x\rA_{L^\infty}$, such that
\be\label{last1}
\begin{aligned}
\intL \chi \psi_x^2 \dx &\le 
A\intL \chi (G(\eta)\psi)^2\dx +A\intL \phi_x^2\arrowvert_{y=-h}\dx\\
&\quad-A\intL\int_{-h}^{\eta(t,x)}\chi_x\phi_x\phi_y\dydx.
\end{aligned}
\ee
Notice that this result holds in fact for any smooth cut-off function $\chi$ and one can also 
take $\chi=1$. 
In this case this inequality simplifies since the last term in the right-hand side vanishes. This gives 
a way to control the $L^2$-norm of $\psi_x$ by the $L^2$-norm of $G(\eta)\psi$. When $\eta$ is smooth, 
this can be obtained by delicate commutator estimates. Here we will give a simple proof which applies 
for any Lipschitz domain.

\section{An integral inequality}

As explained in Section~\ref{S:23}, the first key step is to obtain an inequality 
which involves the integral in time of the energy. 

\begin{prop}\label{T1bis}
Let $m\in C^\infty(\xT)$ be a smooth $2L$-periodic function which is odd and such that $m(L)=0$ (in particular $m(-L)=m(0)=0$). 
Consider a regular solution of the water-wave system defined on the time interval $[0,T]$ 
(see Definition~\ref{D1}) and set 
$$
\zeta=\partial_x(m\eta)+\tdm (1-m_x)\eta-\uq \eta,\quad 
\rho=\zeta + \eta-x\eta_x.
$$
Assume that, for all time $t\in [0,T]$ and all $x\in [-L,L]$, the following assumptions hold:
\begin{alignat*}{2}
&i)\qquad &&\rho(t,x) \ge -\frac{h}{4},\qquad  \la \rho_x(t,x)\ra< \frac 1 4,\\
&ii) \qquad &&\int_{-L}^L (1-m_x(x))\eta(t,x) \dx \le \frac{hL}{3},\qquad 
\la m_x(x)\ra\la \eta_x(t,x)\ra^2\le 2,\\
&iii)\qquad &&\kappa m_{xx}(x)^2\le g, \qquad m_x(x)\le 1.
\end{alignat*}
Then,
\be\label{t49}
\begin{aligned}
&\frac 1 4 \int_0^T \Hr(t)\dt+\frac{h}{4}\int_0^T\int_{-L}^L\phi_x(t,x,-h)^2\dxdt\\
&\qquad\qquad\qquad \le 
\int_0^T\int_{-L}^L \left(\tdm(1-m_x)\widetilde{\psi}+(x-m)\psi_x\right)G(\eta)\psi\dxdt \\
&\qquad\qquad\qquad \quad -\int_0^T\int_{-L}^L P_{ext}\, \zeta\, \dxdt-\int_{-L}^L \zeta \widetilde{\psi} \dx\ST,
\end{aligned}
\ee
with
\be\label{t62}
\widetilde{\psi}(t,x)=\psi(t,x)-\frac{1}{2L}\int_{-L}^L \psi(t,x)\dx,
\ee
where $\int_{-L}^L f\dx \ST$ is a shorthand notation for $\int_{-L}^L f(T,x)\dx-\int_{-L}^L f(0,x)\dx$.
\end{prop}
\begin{rema}
This result holds for any regular solutions (that is, without requiring that $P_{ext}$ is related to $(\eta,\psi)$). 
Also, the assumptions on $m$ are less restrictive then the ones given by Definition~\ref{defi:P}.
\end{rema}
\begin{proof}
The proof follows the analysis in \cite{Boundary,A-Stab-WW}. The 
main novelty is that we are now able to take into account surface tension. 

\begin{nota*}
We write simply
$$
\int \, dx ,\quad \int \, dy ,\quad  
\int\,dt,
$$
as shorthand notations for, respectively,
$$
\int_{-L}^L \, dx ,\quad \int_{-h}^{\eta(t,x)} \, dy ,\quad 
\int_0^T \,dt.
$$
Then, $\iint f\dxdt =\int_0^T\int_{-L}^L f(t,x)\dxdt$ together with similar conventions for other 
integrals.
\end{nota*}

The proof is in four steps. The first step is to use the multiplier method to obtain an 
identity.  We use the following variant of the multiplier method (as 
introduced in our previous paper \cite{Boundary}): we set
\be\label{I22}
A\defn \iint \big\{(\partial_t\eta)(m\px\psi)-(\partial_t\psi)(m\px\eta)\big\}\,\dxdt,
\ee
and compute $A$ in two different ways.

\begin{lemm}\label{T1}
There holds,
\begin{align*}
&\iint m_x \left(- \eta\partial_t\psi-\frac{g}{2}\eta^2-\kappa \frac{1}{\sqrt{1+\eta_x^2}}\right) \,\dxdt+\int F(t)\dt\\
&\qquad\qquad =-\int \px(m \eta) \psi\, \dx\ST-\iint P_{ext}m\eta_x\, \dxdt,
\end{align*}
where
\be\label{t31b}
F(t)=\int (G(\eta)\psi)m\psi_x\dx+ \int (N(\eta)\psi)m\eta_x\dx.
\ee
\end{lemm}
\begin{proof}
Since $m(-L)=m(L)=0$, 
integrating by parts in space and time, we find that $A$ (as given by \e{I22}) satisfies
\be\label{t32}
A=-\int \px(m \eta) \psi\dx\ST+\iint m_x \eta\partial_t\psi \,\dxdt.
\ee
Now, we compute $A$ by replacing $\partial_t\eta$ and $\partial_t\psi$ by the expressions given by System~\e{systemT}. We find that
\be\label{t34}
A=\iint \left(P_{ext}+g\eta-\kappa \px\left(\frac{\eta_x}{\sqrt{1+\eta_x^2}}\right)\right)m\eta_x\dxdt+\int F(t)\dt,
\ee
where $F$ is given by \e{t31b}. Since $m(-L)=m(L)=0$, integrating by parts, we obtain
$$
-\int g\eta m\eta_x\dx=\mez \int g m_x \eta^2\dx.
$$
Similarly, integrating by parts in $x$, one can write
\begin{align*}
-\iint \px\left(\frac{\eta_x}{\sqrt{1+\eta_x^2}}\right)m\eta_x\dx&=\int m_x \frac{\eta_x^2}{\sqrt{1+\eta_x^2}}\dx 
+\int \frac{\eta_x \eta_{xx}}{\sqrt{1+\eta_x^2}} m \dx\\
&=\int m_x \frac{\eta_x^2}{\ai} +\int m\px \ai \dx\\
&=\int m_x \left(\frac{\eta_x^2}{\ai}-\ai\right)\dx\\
&=-\int m_x \frac{1}{\ai}\dx.
\end{align*}
By combining this identity with \e{t32} and \e{t34}, we obtain the wanted result.
\end{proof}

Our next step is to transform the previous identity to make appear a quantity which controls 
the potential energy. When there is no surface tension, it is easy to do so. 
By contrast, with surface tension, we obtain an identity with a remainder (see $R_1$ below) 
which will be estimated later on.

\begin{lemm}\label{T2}
There holds
\be\label{C1}
\begin{aligned}
&\iint \left(\frac g 2 \eta^2+\tdm \kappa \frac{\eta_x^2}{\ai}\right)\dxdt+\int R_1(t)\dt\\
&\qquad\qquad=\iint (1-m_x)\left(\frac g 2 \eta^2+\tdm \kappa \frac{\eta_x^2}{\ai}\right)\dxdt\\
&\qquad\qquad\quad-\int\px(m \eta) \psi\, \dx\ST-\iint P_{ext}\px(m\eta)\, \dxdt\\
&\qquad\qquad\quad-\iint (G(\eta)\psi)(m\psi_x)\dxdt -\iint (N(\eta)\psi) \px (m\eta)\dxdt,
\end{aligned}
\ee
where $R_1$ is given by 
\be\label{tt31b}
R_1(t)=\kappa \int m_{xx} \frac{\eta\eta_x}{\ai}\dx
+\kappa\int m_x \left(1-\mez \frac{\eta_x^2}{\ai}-\frac{1}{\ai}\right)\dx.
\ee
\end{lemm}
\begin{proof}
Directly from the equation for $\psi$, we have
\begin{multline*}
\iint -m_x \eta \partial_t \psi\dxdt=\\
=\iint m_x\left(g\eta^2 +\eta N(\eta)\psi-\kappa \eta \curv  +\eta P_{ext}\right)\dxdt.
\end{multline*}
Noticing that $\eta_x$ vanishes for $x=-L$ or $x=L$ (since $\eta$ is $C^1$, even and $2L$-periodic), we can integrate 
by parts in the term involving the surface tension, to obtain
\begin{align*}
\iint -m_x \eta \partial_t \psi\dxdt&=\iint \left(gm_x \eta^2 
+m_x \eta N(\eta)\psi
+m_x \eta P_{ext}\right)\dxdt\\
&\quad +\iint \kappa m_x \frac{\eta_x^2}{\ai}\dxdt+\iint \kappa m_{xx}\frac{\eta \eta_x}{\ai}\dxdt.
\end{align*}
By combining this identity with Lemma~\ref{T1}, we deduce that
\be\label{a1}
\begin{aligned}
&\iint m_x \left( \frac{g}{2}\eta^2 +\kappa \frac{\eta_x^2}{\ai}-\kappa \frac{1}{\ai}\right)\dxdt
+\kappa \iint m_{xx}\frac{\eta \eta_x}{\ai}\dxdt\\
&\qquad\qquad=-\iint P_{ext} \px(m\eta) \dxdt -\int \px (m\eta)\psi \dx \ST\\
&\qquad\qquad\quad -\iint (G(\eta)\psi)(m\psi_x)\dxdt -\iint (N(\eta)\psi) \px (m\eta)\dxdt.
\end{aligned}
\ee
Since $m(-L)=m(L)=0$, we have $\int m_x \dx=0$ which implies that
\begin{multline*}
-\iint m_x \frac{1}{\ai}\dxdt=\\
=\iint m_x\left( \mez \frac{\eta_x^2}{\ai}+\left[ -\mez \frac{\eta_x^2}{\ai}-\frac{1}{\ai}+1\right]\right)\dxdt.
\end{multline*}
By reporting this identity in \e{a1}, we conclude that
\be\label{C20}
\begin{aligned}
&\iint m_x \left(\frac{g}{2}\eta^2+\frac{3\kappa}{2}\frac{\eta_x^2}{\ai}\right)\dxdt+\int R_1\dt=\\
&\qquad\qquad=-\iint P_{ext} \px(m\eta) \dxdt -\int \px (m\eta)\psi \dx \ST\\
&\qquad\qquad\quad -\iint (G(\eta)\psi)(m\psi_x)\dxdt -\iint (N(\eta)\psi) \px (m\eta)\dxdt.
\end{aligned}
\ee
We then deduce the wanted result by splitting $m_x$ as $1+(m_x-1)$ in the left-hand side.
\end{proof}

The previous lemma gives an identity whose left hand-side is related to the integral in time of the 
potential energy. Since in the end we need to control the energy (and not only the potential energy), 
we need to compare this term to a similar one related to the kinetic energy. 
This is the purpose of the following lemma which, loosely speaking, 
is a version of the principle of equipartition of energy. 
Let us recall that the kinetic energy is given by
$$
\mez\int_{-L}^L\int_{-h}^{\eta(t,x)}\la \nabla_{x,y}\phi \ra^2\dydx=\mez\intL \psi G(\eta)\psi \dx,
$$
where we used the divergence theorem and the fact that $\phi$ is harmonic.

\begin{lemm}\label{L3}
For any smooth function $\theta=\theta(x)$, there holds
\be\label{t70rho}
\begin{aligned}
\iint \theta \left(g\eta^2+\kappa\frac{\eta_x^2}{\ai}\right)\dxdt
&= \iint \theta\psi G(\eta)\psi\dxdt\\
&\quad- \iint \theta \eta P_{ext}\dxdt 
- \int\theta \eta\psi\dx\ST\\
&\quad -\iint \theta \eta \, N(\eta)\psi\dxdt\\
&\quad -\iint \kappa \theta_x \frac{\eta \eta_x}{\ai}\dxdt.
\end{aligned}
\ee
\end{lemm}
\begin{proof}
We have
$$
\iint \theta\psi G(\eta)\psi\dxdt=\iint \theta \psi \partial_t \eta \dxdt
=\int \theta\psi\eta\dx\ST-\iint \theta\eta\partial_t\psi \dxdt.
$$
So \e{t70rho} easily follows from the equation for $\psi$ (integrating by parts as above in 
the term involving the surface tension).
\end{proof}

Recall that the energy is given at time $t$ by
$$
\mathcal{H}(t)=\frac{g}{2}\int\eta^2\dx+\kappa \int \left(\sqrt{1+\eta_x^2}-1\right)\dx
+\mez \int \psi G(\eta)\psi\dx.
$$
Instead of bounding $\int \Hr \dt$, we will estimate the integral in time of the following quantity: 
$$
\widetilde{H}(t)\defn\int \left( \frac g 2 \eta^2+\kappa \frac{\eta_x^2}{\ai}+\mez \psi G(\eta)\psi\right)\dx.
$$
It will be sufficient to estimate this term since $\widetilde{H}(t)\ge \mathcal{H}(t)$. 
Indeed, 
$$
\sqrt{1+\eta_x^2}-1=\frac{\eta_x^2}{1+\sqrt{1+\eta_x^2}}\le \frac{\eta_x^2}{\ai}.
$$

We now combine the previous results to obtain an identity for the integral in time of $\widetilde{H}(t)$. 

\begin{lemm}\label{L:36}
Set
$$
\zeta=\partial_x(m\eta)+\tdm (1-m_x)\eta-\uq \eta.
$$
Then
\be\label{C2}
\begin{aligned}
\mez \int \widetilde{H}(t)\dt+\int R_2(t)\dt&=
\tdm \iint (1-m_x)\psi G(\eta)\psi \dxdt\\
&\quad -\iint P_{ext}\zeta \dxdt-\int \psi \zeta\dx \ST\\
&\quad-\iint m\psi_x G(\eta)\psi \dxdt\\
&\quad -\iint \zeta N(\eta)\psi\dxdt,
\end{aligned}
\ee
where 
\begin{align*}
R_2(t)&=-\frac \kappa 2 \int m_{xx}\frac{\eta \eta_x}{\ai}\dx \\
&\quad +\kappa\int m_x\left(1-\frac{\eta_x^2}{2\ai}
-\frac{1}{\ai}\right)\dx\\
&\quad +\frac{3\kappa}{4} \int \frac{\eta_x^2}{\ai}\dx +g\int (1-m_x)\eta^2\dx.
\end{align*}
\end{lemm}
\begin{proof}
We will deduce this identity from \e{C1} and from Lemma~\ref{L3}. 
Denote by $I_1$ the first term of the left-hand side in \e{C1}, given by
$$
I_1=\iint \left(\frac g 2 \eta^2+\tdm \kappa \frac{\eta_x^2}{\ai}\right)\dxdt.
$$
We also introduce 
$$
I_m=\iint (1-m_x)\left(\frac g 2 \eta^2+\tdm \kappa \frac{\eta_x^2}{\ai}\right)\dxdt.
$$
With these notations, the identity \e{C1} reads
\be\label{C3}
\begin{aligned}
I_1+\int R_1(t)\dt&=I_m-\intL\px(m \eta) \psi\, \dx\ST-\iint P_{ext}\px(m\eta)\, \dxdt\\
&\quad -\iint m\psi_x G(\eta)\psi \dxdt-\iint \px(m\eta) N(\eta)\psi\dxdt.
\end{aligned}
\ee

We split these terms as $I_1=J_1+P_1$ and $I_m=J_m+P_m$ with
\begin{align*}
J_1&=\mez \iint \left(g  \eta^2+ \kappa \frac{\eta_x^2}{\ai}\right)\dxdt+\frac{\kappa}{4}\iint \frac{\eta_x^2}{\ai}\dxdt,\\
P_1&=\frac{3\kappa}{4}\iint \frac{\eta_x^2}{\ai}\dxdt,\\
J_m&=\tdm \iint (1-m_x)\left(g \eta^2+\kappa \frac{\eta_x^2}{\ai}\right)\dxdt,\\ 
P_m&=-\iint (1-m_x)g\eta^2\dxdt.
\end{align*}
We now use Lemma~\ref{L3} to compute the first term contributing to $J_1$. 
This gives
\begin{align*}
\iint \left(g  \eta^2+ \kappa \frac{\eta_x^2}{\ai}\right)\dxdt&=\iint\psi G(\eta)\psi-\int \psi\eta\dx\ST\\
&\quad-\iint \eta N(\eta)\psi \dxdt-\iint \eta P_{ext}\dxdt.
\end{align*}
Consequently,
\begin{align*}
J_1=&\mez \iint \left(\frac g 2 \eta^2+\mez \kappa \frac{\eta_x^2}{\ai}\right)\dxdt+\frac{\kappa}{4}\iint \frac{\eta_x^2}{\ai}\dxdt
\\
&\quad+\uq \iint \psi G(\eta)\psi \dxdt\\
&\quad -\uq \iint \eta P_{ext}\dxdt-\uq \iint \eta N(\eta)\psi\dxdt-\uq \int \psi\eta\dx\ST.
\end{align*}
Then, one checks that
$$
J_1=\mez \int \widetilde{H}\dt -\uq \iint \eta P_{ext}\dxdt-\uq \iint \eta N(\eta)\psi\dxdt-\uq \int \psi\eta\dx\ST.
$$
We can also use Lemma~\ref{L3} to obtain that
\begin{align*}
J_m=&\tdm \iint (1-m_x)\psi G(\eta)\psi \dxdt\\
&\quad -\tdm \iint (1-m_x)\eta P_{ext}\dxdt-\tdm \iint (1-m_x)\eta N(\eta)\psi\dxdt\\
&\quad-\tdm \int (1-m_x)\psi\eta\dx\ST\\
&\quad+\tdm \kappa \iint m_{xx}\frac{\eta \eta_x}{\ai}\dxdt.
\end{align*}
By combining the previous identities with \e{C3} we obtain the wanted result \e{C2}. 
\end{proof}

Our next task is to study the 
last two terms in the right-hand side of \e{C2},
$$
\int m\psi_x G(\eta)\psi \dx,\quad 
\int \zeta N(\eta)\psi\dx.
$$
The following lemma, which relies crucially on the results proved in \cite{Boundary}, 
shows that these two terms can either be recast in a simpler way, or give rise to signed quantities (with the good signs).  

\begin{lemm}
Set
$$
\rho=\zeta+\eta-x\eta_x.
$$
There holds
\be\label{C5}
\begin{aligned}
\mez \int \widetilde{H}(t)\dt+\int R_3(t)\dt&=
\tdm \iint (1-m_x)\psi G(\eta)\psi \dxdt\\
&\quad -\iint P_{ext}\zeta \dxdt-\int \psi \zeta\dx \ST\\
&\quad-\iint (m-x)\psi_x G(\eta)\psi \dxdt\\
&\quad +\iiint \rho_x \phi_x\phi_y \dydxdt,
\end{aligned}
\ee
where 
$$
R_3(t)=R_2(t)+\mez \int_{-L}^L (h+\rho) \phi_x^2\arrowvert_{y=-h}\dx +L\int_{-h}^{\eta(t,L)}   \phi_y^2\arrowvert_{x=L}\dy.
$$
\end{lemm}
\begin{proof}
To handle the term $\int m\psi_x G(\eta)\psi \dx$ 
we split $m\psi_x$ as $m\psi_x=x\psi_x+(m-x)\psi_x$ to obtain an expression which 
makes appear a positive term through a Pohozaev identity. We will use the following identity 
(proved in \cite{Boundary})
$$
\int_0^L (G(\eta)\psi)x\psi_x\dx=\widetilde\Sigma+ \int_0^L (\eta-x\eta_x)N(\eta)\psi\dx,
$$
where $\widetilde\Sigma(t)$ is a positive term given by
$$
\widetilde\Sigma(t)=\frac{h}{2} \int_0^L \phi_x^2(t,x,-h)\dx+\frac{L}{2}\int_{-h}^{\eta(t,L)}\phi_y^2(t,L,y)\dy. 
$$
Owing to the fact that $\eta$, $\psi$, $\phi_x^2$ and $\phi_y^2$ are even in $x$, we deduce that
\be\label{C4}
\int_{-L}^L (G(\eta)\psi)x\psi_x\dx=\Sigma+ \int_{-L}^L (\eta-x\eta_x)\big(N(\eta)\psi\big)\dx,
\ee
where
$$
\Sigma(t)=\frac{h}{2} \int_{-L}^L \phi_x^2(t,x,-h)\dx+L\int_{-h}^{\eta(t,L)}\phi_y^2(t,L,y)\dy. 
$$
It follows from \e{C4} that
$$
\intL m\psi_x G(\eta)\psi \dx=\Sigma(t)+\intL (\eta-x\eta_x) N(\eta)\psi\dx+\intL (m-x)\psi_x G(\eta)\psi\dx.
$$

Recall that $\rho=\zeta+\eta-x\eta_x$. 
To complete the proof it remains only to show that one can write the integral $\intL \rho N(\eta)\psi\dx$ 
as the sum of two terms which can be controlled by the sum of the energy and 
the positive term $\Sigma$ given by the Pohozaev identity. 
To do so, we use the following identity (which is proved in the appendix, cf \e{C6-bis}), 
\be\label{C6}
\intL \rho N(\eta)\psi\dx
=-\intL\int_{-h}^{\eta(t,x)} \rho_x\phi_x\phi_y \dydx +\mez \intL \rho \phi_x^2\arrowvert_{y=-h}\dx.
\ee
By plugging this result into \e{C2}, we complete the proof of the lemma.
\end{proof}

Recall that, by notation,
$$
\widetilde\psi(t,x)=\psi(t,x)-\langle\psi\rangle(t)\quad \text{with}\quad \langle\psi\rangle(t)=\frac{1}{2L}\int_{-L}^L \psi(t,x)\dx.
$$
We now have to prove that one can replace $\psi$ by $\widetilde{\psi}$ in the previous identity 
\e{C5}, to the price of an admissible error term. 

\begin{lemm}\label{L2.7}
There holds
\be\label{C7}
\begin{aligned}
&\mez \int \widetilde{H}(t)\dt+\int R_4(t)\dt\\
&\qquad\qquad=
 \iint \left(\tdm(1-m_x)\widetilde{\psi} G(\eta)\psi +(x-m)\psi_x G(\eta)\psi\right)\dxdt\\
&\qquad\qquad\quad -\iint P_{ext}\zeta \dxdt-\int \widetilde{\psi} \zeta\dx \ST\\
&\qquad\qquad\quad +\iiint \rho_x \phi_x\phi_y \dydxdt,
\end{aligned}
\ee
where 
$$
R_4(t)=R_3(t)-\beta(t)\int \phi_x^2\arrowvert_{y=-h}\dx\quad\text{with}\quad 
\beta(t)=\frac{3}{8L} \int (1-m_x)\eta \dx.
$$
\end{lemm}
\begin{proof}
Introduce 
\be\label{defiI}
I=\iint \tdm (1-m_x)\langle\psi\rangle G(\eta)\psi\dxdt-\int \langle\psi\rangle\zeta\dx\ST.
\ee
Then
$$
\int R_3(t)\dt-I=\int R_4(t)\dt.
$$
Since $\int \eta\dx=0$, by definition of $\zeta$, we have
$$
\int \zeta \dx=\tdm \int (1-m_x)\eta\dx.
$$
Since $\partial_t\eta=G(\eta)\psi$, we deduce from \e{defiI} that
$$
I=\iint \tdm (1-m_x)\langle\psi\rangle\partial_t \eta\dxdt-\int \tdm(1-m_x)\eta \langle\psi\rangle\dx\ST.
$$
Integrating by parts in time, this implies that
$$
I=-\tdm \iint (1-m_x)\eta\partial_t \langle\psi\rangle\dxdt,
$$
and hence to compute $I$ we need only to compute $\partial_t \langle\psi\rangle$. To do so, recall that
$$
\partial_t \psi+g\eta+N(\eta)\psi-\kappa \px \left(\frac{\eta_x}{\ai}\right)+P_{ext}=0.
$$
Recall that $\int P_{ext}\dx=0$ and $\int \eta\dx=0$ by assumptions 
(see \e{t9b}). On the other hand, since $\eta$ is $C^1$ and even in $x$, one has 
$\eta_x(t,-L)=\eta_x(t,L)$ and hence we conclude that $\langle\psi\rangle=(2L)^{-1}\int_{-L}^L \psi\dx$ satisfies
$$
\partial_t \langle\psi\rangle=-\frac{1}{2L}\int N(\eta)\psi\dx.
$$
Now we use that the identity \e{C6-bis} with $\mu=1$, to obtain that
$$
\partial_t \langle\psi\rangle=-\frac{1}{4L}\int \phi_x^2\arrowvert_{y=-h}\dx.
$$
As a result,
$$
I=\int \beta(t)\int \phi_x^2\arrowvert_{y=-h}\dxdt \quad\text{with}\quad 
\beta(t)=\frac{3}{8L} \int (1-m_x)\eta \dx.
$$
This completes the proof of the lemma.
\end{proof}

To complete the proof of the proposition, it remains to estimate the remainder terms. 
This is the purpose of the following lemma.

\begin{lemm}\label{CL8}
Assume that, for all time $t\in [0,T]$ and all $x\in [0,L]$, the following assumptions hold:
\begin{alignat*}{2}
&i)\qquad &&\rho(t,x) \ge -\frac{h}{4},\qquad  \la \rho_x(t,x)\ra\le \frac 1 4,\\
&ii) \qquad &&\int_{-L}^L (1-m_x(x))\eta(t,x) \dx \le \frac{hL}{3},\qquad 
\la m_x(x)\ra\la \eta_x(t,x)\ra^2\le 2,\\
&iii)\qquad &&\kappa m_{xx}(x)^2\le g,\qquad m_x(x)\le 1.
\end{alignat*}
Then, for all $t\in [0,T]$,
$$
R_4(t)-\iint_{\Omega(t)} \rho_x \phi_x\phi_y \dydx\ge \frac{h}{4}\int\phi_x^2\arrowvert_{y=-h}\dx-\uq \widetilde{H}(t).
$$
\end{lemm}
\begin{proof}
We have $R_4=I_1+\cdots+I_6$ with
\begin{alignat*}{2}
I_1&=\int \left(\frac{h+\rho}{2}-\beta\right)\phi_x^2\arrowvert_{y=-h}\dx,\quad
&& I_2=L\int \phi_y^2\arrowvert_{x=L}\dy,\\
I_3&=\frac{3\kappa}{4} \int\frac{\eta_x^2}{\ai}\dx,\qquad &&I_4=g\int (1-m_x)\eta^2\dx,\\
I_5&=-\frac{\kappa}{2}\int m_{xx}\frac{\eta \eta_x}{\ai}\dx,\quad && \\
I_6&=\kappa \int m_x\left(1-\frac{\eta_x^2}{2\ai}-\frac{1}{\ai}\right)\dx &&.
\end{alignat*}
By assumption $i)$ and $ii)$ we have $\rho\ge -h/4$ and $\beta\le h/8$, so
$$
I_1\ge \frac{h}{4}\int\phi_x^2\arrowvert_{y=-h}\dx.
$$ 
Obviously we have $I_2\ge 0$ and $I_4\ge 0$ since $m_x(x)\le 1$ by assumption $iii)$. 
So, to prove the lemma, we need only to show that
\be\label{C10}
I_3+I_5+I_6-\iint_{\Omega(t)} \rho_x \phi_x\phi_y \dydx\ge -\uq \widetilde{H}(t).
\ee
We begin by estimating $I_6$. Introduce $f\colon \xR_+\rightarrow \xR$ defined by
$$
f(u)=1-\frac{u}{2\sqrt{1+u}}-\frac{1}{\sqrt{1+u}}.
$$
Since $f(0)=0$, we have the estimate
$$
\la f(u)\ra\le \int_0^u \la f'(s)\ra \, {\rm d}s\le \uq \int_0^u s \, {\rm d}s=\frac{u^2}{8}.
$$
So, since $|\eta_x|\le 1$ by assumption, we may write
$$
\la I_6\ra\le \kappa \int \la m_x\ra \frac{\eta_x^4}{8}\dx \le \frac{\kappa}{4}\int \la m_x\ra\frac{\eta_x^2}{\sqrt{1+\eta_x^2}}\eta_x^2\dx.
$$
Now we use the assumption $\la m_x(x)\ra\la \eta_x(t,x)\ra^2\le 2$ to infer that
$$
\la I_6\ra\le \frac{\kappa}{2}\int \frac{\eta_x^2}{\sqrt{1+\eta_x^2}}\dx.
$$
We now move to $I_5$. Writing
\begin{align*}
\la I_5\ra&\le \frac{\kappa}{2}\int \uq m_{xx}^2\frac{\eta^2}{\ai}\dx +\frac{\kappa}{2}\int \frac{\eta_x^2}{\ai}\dx\\
&\le \frac{\kappa}{8}\sup_{[-L,L]}\la m_{xx}\ra^2 \int \eta^2\dx +\frac{\kappa}{2}\int \frac{\eta_x^2}{\ai}\dx,
\end{align*}
and using assumption $iii)$, one obtains that
$$
\la I_5\ra\le \frac{g}{8}\int \eta^2\dx+\frac{\kappa}{2}\int \frac{\eta_x^2}{\ai}\dx.
$$
On the other hand,
\begin{align*}
\iint_{\Omega(t)} \rho_x \phi_x\phi_y \dydx&\le 
\mez \sup_{[-L,L]}\la \rho_x\ra \iint_{\Omega(t)}(\phi_x^2+\phi_y^2)\dydx\\
&\le \frac{1}{8} \iint_{\Omega(t)}\la \nabla_{x,y}\phi\ra^2\dydx,
\end{align*}
since $\la \rho_x\ra\le 1/4$ by assumption. 
By combining the previous inequalities with the definition of $I_3$, we conclude that the left-hand side of \e{C10} is greater than
$$
-\frac{g}{8}\int \eta^2\dx-\frac{\kappa}{4}\int \frac{\eta_x^2}{\ai}\dx -\frac{1}{8} \iint_{\Omega(t)}\la \nabla_{x,y}\phi\ra^2\dydx
=-\uq \widetilde{H}(t).
$$
This completes the proof of the lemma.
\end{proof}

In view of Lemma~\ref{L2.7} and Lemma~\ref{CL8}, we obtain 
Proposition~\ref{T1bis}, recalling that $\widetilde{H}(t)\ge \mathcal{H}(t)$ 
as explained above the statement of Lemma~\ref{L:36}. 
\end{proof}

\section{Stabilization}

In this section we prove Theorem~\ref{T:2}. 
In the previous section, we proved an inequality which is valid 
for any value of $P_{ext}$. Here, by contrast, we will prove an inequality which 
exploits in a crucial way some special choice for $P_{ext}$. 
We now assume that the external pressure is given by
$$
P_{ext}(t,x)=P(t,x) +p(t),
$$
where
$$
P(t,x)=\lambda\chi(x)\partial_t\eta,\quad 
p(t)=-\frac{1}{2L}\int_{-L}^L P(t,x)\dx,
$$
for some positive constant $\lambda$ and some cut-off function $\chi$ to be chosen. 
In the introduction, to state Theorem~\ref{T:2}, we considered only the case $\lambda=1$ (to fix ideas) and indeed 
we will not try to optimize $\lambda$. However, in this section 
we add this parameter for possible applications where it 
could be important to tune it. 
 
Notice that, since $\partial_t\eta=G(\eta)\psi$, we have
$$
P(t,x)=\lambda\chi(x)G(\eta)\psi.
$$
Let us recall also how $\chi$ is defined. 
\begin{defi}
Fix $\delta>0$ and consider a $2L$-periodic $C^\infty$ function $\cutoff$ satisfying $0\le \cutoff\le 1$ and such that
$$
\cutoff(x)=\cutoff(-x),\quad x\cutoff'(x)\le 0 \text{ for } x\in [-L,L],\quad
\cutoff(x)=\left\{\begin{aligned}&1 \text{ if }x\in [0,L-\delta],\\ &0 \text{ if }x\in \left[L-\frac{\delta}{2},L\right].\end{aligned}\right.
$$
We set  
$$
m(x)=x\cutoff(x) \quad \text{and}\quad \chi(x)=1-m_x(x).
$$
\end{defi}

As already mentioned, it follows from the hamiltonian structure of the equation that
$$
\fract \Hr(t)=-\int_{-L}^L P_{ext} G(\eta)\psi\dx.
$$
In addition, since $\int_{-L}^L G(\eta)\psi\dx=0$, we get $\int_{-L}^L p(t)G(\eta)\psi\dx=0$ and hence
\be\label{C14}
\fract \Hr(t)=-\int_{-L}^L P G(\eta)\psi\dx=-\lambda\int_{-L}^L \chi(G(\eta)\psi)^2\dx.
\ee
As a result, since $\lambda>0$ and $\chi\ge 0$, we verify that the energy $\Hr$ is a decreasing function of time. 
Our goal in this section is to obtain a quantitative result which asserts that, 
if the mild assumptions on $\eta$ given by the statement of Theorem \ref{T:2} are satisfied, 
then there exists a positive constant $C$, depending only on $\kappa,g,h,L$, such that
\be\label{d1}
\Hr(T)\le \frac{C}{T} \Hr(0).
\ee
To do so, as already explained in the introduction, we will prove that 
\be\label{C11}
\int_0^T\Hr(t)\dt\le C\Hr(0).
\ee
This will imply the desired result \e{d1}, by writing
$$
\Hr(T)\le \frac{1}{T} \int_0^T\Hr(t)\dt\le \frac{C}{T}\Hr(0).
$$

We begin by establishing two properties of the pressure which explain the choices of the expressions for 
$P$, $P_{ext}$ and $\chi$.
\begin{lemm}
There holds
\be\label{d7}
\intTL P^2\dxdt \le \lambda (\sup_{[-L,L]} \chi) \Hr(0),
\ee
and
\be\label{d8}
-\intT p(t) \intL \zeta(t,x)\dxdt\le \frac{3\lambda}{2g} \sup_{[-L,L]}(1-m_x)^2 \Hr(0),
\ee
where recall that $\zeta$ is given by
$$
\zeta=\partial_x(m\eta)+\tdm (1-m_x)\eta-\uq \eta.
$$
\end{lemm}
\begin{proof}
Firstly, integrating in time the identity \e{C14}, we obtain 
\be\label{Bures4}
\lambda \intTL \chi(G(\eta)\psi)^2\dxdt =\Hr(0)-\Hr(T) \le  \Hr(0).
\ee
This immediately implies that
$$
\intTL P^2\dxdt=\lambda^2 \intTL \chi^2 \big(G(\eta)\psi\big)^2\dxdt \le \lambda(\sup \chi)\Hr(0),
$$
which is \e{d7}. 

To estimate $\iint p\zeta \dxdt$ we use in a crucial way the fact that, by definition, $\chi(x)=1-m_x(x)$. 
Since $\int \eta\dx=0$, directly from 
the definition of $\zeta$, we have
$$
\intL \zeta(t,x)\dx=\tdm\nu(t)\quad \text{where}\quad \nu(t)\defn \intL (1-m_x)\eta\dx. 
$$
So 
$$
\intL p(t) \zeta(t,x) \dx = \tdm p(t)\nu(t).
$$
On the other hand, by definition,
$$
p=-\frac{\lambda}{2L}\intL \chi G(\eta)\psi\dx=-\frac{\lambda}{2L}\intL \chi \partial_t \eta\dx =-\frac{\lambda}{2L}\partial_t \intL \chi \eta\dx.
$$
Since $\chi=1-m_x$, we deduce that
$$
p=-\frac{\lambda}{2L}\partial_t \nu,
$$
and hence
$$
\intT p(t) \intL \zeta(t,x)\dxdt=-\frac{3\lambda}{4L} \intT \nu \partial_t \nu \dt.
$$
Consequently,
$$
-\intTL p(t)\zeta(t,x)\dxdt=\frac{3\lambda}{8 L}\int_0^T \fract\nu^2\dt \le\frac{3\lambda}{8 L}\nu(T)^2.
$$
Now, to estimate $\nu(T)^2$ we use the Cauchy-Schwarz inequality to obtain
\begin{align*}
\nu(T)^2&=\left(\int_{-L}^L (1-m_x)\eta\dx\right)^2\le 
2L \int_{-L}^L (1-m_x)^2\eta^2\dx\\
&\le 2L \sup_{[-L,L]}(1-m_x)^2 \int_{-L}^L \eta^2\dx\\
&\le \frac{4L}{g} \sup_{[-L,L]}(1-m_x)^2 \Hr(T)\le \frac{4L}{g} \sup_{[-L,L]}(1-m_x)^2 \Hr(0),
\end{align*}
since $\Hr$ is decreasing. 
This yields 
$$
-\intTL p(t)\zeta(t,x)\dxdt\le \frac{3\lambda}{2g}\sup_{[-L,L]}(1-m_x)^2\, \Hr(0),
$$
which completes the proof of the lemma.
\end{proof}

We now observe that the assumptions of Theorem~\ref{T:2} are stronger than the ones of Proposition~\ref{T1bis} so we may apply this proposition 
to write that
\be\label{C12}
\frac 1 4 \int_0^T \Hr(t)\dt\le O+W+B-I,
\ee
where 
\begin{align*}
I&\defn \frac{h}{4}\int_0^T\int_{-L}^L\phi_x(t,x,-h)^2\dxdt,\\
O&\defn \int_0^T\int_{-L}^L \left(\tdm(1-m_x)\widetilde{\psi}+(x-m)\psi_x\right)G(\eta)\psi\dxdt, \\
W&\defn -\int_0^T\int_{-L}^L P_{ext}\, \zeta\, \dxdt,\\
B&\defn -\int_{-L}^L \zeta \widetilde{\psi} \dx\ST.
\end{align*}
To deduce \e{C11} from \e{C12}, as already explained, it is sufficient to prove that that there exists a constant $K$ depending only on $g,\kappa,h,L$ such that 
\be\label{C13}
O+W+B -I \le K \Hr(0)+a \int_0^T\Hr(t)\dt \quad\text{for some}\quad a <\uq.
\ee
Indeed, \e{C13} implies
$$
\int_0^T\Hr(t)\dt\le \frac{K}{1/4-a}\Hr(0),
$$
which is the wanted result~\e{C11}.

To prove \e{C13}, we begin by estimating the term $W$.

\begin{lemm}\label{CL7}
Let $\eps_1>0$ and set
\begin{align}
C_1&=\frac{\lambda}{2\eps_1}\sup_{[-L,L]} \chi(x) +\frac{3\lambda}{2g}\sup_{[-L,L]}(1-m_x(x))^2,\\
K_1&=\frac{6}{\kappa}\sup_{[-L,L]} m(x)^2 +\frac{4}{g}\sup_{[-L,L]}\left( \frac{5}{4}-\frac{m_x(x)}{2}\right)^2.\label{Bures1}
\end{align}
Then there holds
\be\label{d9}
-\iint P_{ext} \zeta\dxdt \le C_1\Hr(0)+\frac{\eps_1 K_1}{2}\int \Hr(t)\dt.
\ee
\end{lemm}
\begin{proof}
Recall that $P=P_{ext}+p$. We have already estimated $-\iint p \zeta\dxdt$ so to prove \e{d9} it remains only 
to estimate $-\iint P \zeta\dxdt$. 
Firstly, using \e{d7} and the Young's inequality $-ab\le (1/2\eps_1) a^2+(\eps_1/2)b^2$, 
we have
\begin{align*}
-\iint P \zeta \dxdt &\le \frac{1}{2\eps_1}\iint P^2\dxdt +\frac{\eps_1}{2}\iint \zeta^2\dxdt\\
&\le \frac{\lambda}{2\eps_1}(\sup \chi) \Hr(0)+\frac{\eps_1}{2}\iint \zeta^2\dxdt.
\end{align*}
This implies that $-\iint P \zeta\dxdt$ is bounded by the right-hand side of \e{d9} since 
\be\label{d9bis}
\int \zeta(t,x)^2\dx\le K_1\Hr(t),
\ee
by definition of $\zeta$, $K_1$ and $\Hr$. Indeed,
\begin{align*}
\int \zeta^2\dx &\le \int 2 m^2\eta_x^2\dx +\int 2\left(\frac 5 4 -\frac{m_x}{2}\right)^2\eta^2\dx\\
&\le \frac{6m^2}{\kappa}\int \kappa\left(\sqrt{1+\eta_x^2}-1\right)\dx +\frac{4}{g}\int \left(\frac 5 4 -\frac{m_x}{2}\right)^2 \frac{g}{2}\eta^2\dx,
\end{align*}
where we used that $1+\sqrt{1+\eta_x^2}\le 3$ since $\la\eta_x\ra\le 1$ by assumption.
\end{proof}

We now estimate the term $B$ in \e{C12}.

\begin{lemm}
There exists a constant $K_2$ depending only on $h,L$ 
such that,
\be\label{Bures2}
\intL \widetilde{\psi}(t,x)^2\dx\le K_2 \Hr(t).
\ee
Moreover, 
\be\label{d9b}
\la \int_{-L}^L \zeta \widetilde{\psi} \dx\ST\ra\le (K_1+K_2)\Hr(0),
\ee
where $K_1$ is given by \e{Bures1}.
\end{lemm}
\begin{proof}
The first estimate follows from usual estimates (Poincar\'e's inequality) and we postpone 
its proof to Appendix~\ref{Proof:Bures2}. To prove \e{d9b}, consider $t_0=0$ or $t_0=T$ and write
\begin{align*}
\la \int \zeta(t_0,x) \widetilde{\psi} (t_0,x)\dx\ra
&\le \mez \int \zeta^2(t_0,x)\dx +\mez \int \widetilde{\psi}^2(t_0,x)\dx\\
&\le \mez(K_1+K_2)\Hr(t_0)\le \mez(K_1+K_2)\Hr(0),
\end{align*}
where we used \e{d9bis} and the fact that $\Hr$ is decreasing. This implies \e{d9b}.
\end{proof}

We are now in position to estimate the last integral, which is
$$
O=\int_0^T\int_{-L}^L \left(\tdm(1-m_x)\widetilde{\psi}+(x-m)\psi_x\right)G(\eta)\psi\dxdt.
$$
\begin{lemm}
For any positive real numbers $\eps_2,\eps_3$, there holds
\begin{align*}
\la O\ra&\le \left(\frac{1}{\eps_2\lambda}\sup_{[-L,L]} \chi+\frac{4\eps_3 L}{\lambda}+\frac{L}{2\eps_3\lambda}\right)\Hr(0)
+2\eps_3 L\intL \phi_x^2\arrowvert_{y=-h}\dxdt \\
&\quad + \left(\frac{9\eps_2}{16}K_2 +4\eps_3 L\sup \la \chi_x\ra\right)\int \Hr(t)\dt.
\end{align*}
\end{lemm}
\begin{proof}
We split $O$ as $O_1+O_2$ with
$$
O_1=\iint \tdm(1-m_x)\widetilde{\psi}G(\eta)\psi\dxdt,\quad 
O_2=\iint (x-m)\psi_x G(\eta)\psi\dxdt.
$$
Since $\chi=1-m_x$, for any $\eps_2>0$, one has the estimate
\begin{align*}
\la O_1\ra&\le \frac{9\eps_2}{16}\iint \widetilde{\psi}^2\dxdt+\frac{1}{\eps_2}\iint (1-m_x)^2(G(\eta)\psi)^2\dxdt\\
&\le \frac{9\eps_2}{16}\iint \widetilde{\psi}^2\dxdt+\frac{1}{\eps_2}\sup_{[-L,L]} \chi \iint \chi(G(\eta)\psi)^2\dxdt.
\end{align*}
An application of \e{Bures4} gives the estimate
$$
\la O_1\ra\le \frac{9\eps_2}{16}\iint \widetilde{\psi}^2\dxdt+\frac{1}{\eps_2 \lambda}\big(\sup_{[-L,L]} \chi \big) \mathcal{H}(0).
$$
It follows from \e{Bures2} that
$$
\la O_1\ra\le \frac{9\eps_2}{16}K_2 \int \Hr(t)\dt+\frac{1}{\eps_2 \lambda}\big(\sup_{[-L,L]} \chi \big)\mathcal{H}(0).
$$

The estimate of $O_2$ is more delicate. 
We begin by proving that, for any $x\in [-L,L]$, one has
\be\label{b151}
\la x-m(x)\ra \le L \chi(x).
\ee
To see this, recall that $m$ is as defined in Definition~\ref{defi:P}. 
Since $m$ is odd in $x$, we can assume without loss of generality that $x\in [0,L]$, so that
$$
x-m(x)=x(1-\cutoff(x))\le L(1-\cutoff(x))\le L(1-\cutoff(x)-x\cutoff'(x))=L(1-m_x(x)),
$$
where we used again the definition $\chi(x)=1-m_x(x)$ together with the assumptions $\cutoff(x)\le 1$ and $x\cutoff'(x)\le 0$.

An application of \e{b151} and the Young's inequality implies that, for any $\eps_3>0$,
$$
\la (x-m)\psi_x G(\eta)\psi\ra \le \frac{\eps_3 L}{2}\chi \psi_x^2+\frac{L}{2\eps_3}\chi (G(\eta)\psi)^2. 
$$
Now we claim that if $\la \eta_x\ra\le 1$ then 
\be\label{d11}
\begin{aligned}
\intL \chi \psi_x^2 \dx &\le 
8\intL \chi (G(\eta)\psi)^2\dx +4\intL \phi_x^2\arrowvert_{y=-h}\dx\\
&\quad-8\intL\int_{-h}^{\eta(t,x)}\chi_x\phi_x\phi_y\dydx.
\end{aligned}
\ee
Assume that this is proved. Then an application of the estimate
$$
\iint \chi(G(\eta)\psi)^2\dxdt\le \frac{1}{\lambda}\Hr(0),
$$
will imply that
\begin{align*}	
\la O_2\ra &\le \frac{1}{\lambda}\left(\frac{\eps_3 L}{2}8+\frac{L}{2\eps_3}\right)\Hr(0)+\frac{\eps_3 L}{2}4\intL \phi_x^2\arrowvert_{y=-h}\dxdt \\
&\quad -\frac{\eps_3 L}{2}8\iiint\chi_x\phi_x\phi_y\dydxdt,
\end{align*}
so
\begin{align*}	
\la O_2\ra &\le \left(\frac{4\eps_3 L}{\lambda}+\frac{L}{2\eps_3\lambda}\right)\Hr(0)+2\eps_3 L\intL \phi_x^2\arrowvert_{y=-h}\dxdt \\
&\quad +4\eps_3 L\sup \la \chi_x\ra \iiint \mez \la\nabla_{x,y}\phi\ra^2\dydxdt,
\end{align*}
and hence
\begin{align*}	
\la O_2\ra &\le \left(\frac{4\eps_3 L}{\lambda}+\frac{L}{2\eps_3\lambda}\right)\Hr(0)+2\eps_3 L\intL \phi_x^2\arrowvert_{y=-h}\dxdt \\
&\quad +4\eps_3 L\sup \la \chi_x\ra \int \Hr(t)\dt.
\end{align*}

It remains to prove \e{d11}. 
Introduce the notations
$V=\phi_x\aeta$, $B=\phi_y\aeta$ and recall that, by definition, $G(\eta)\psi=\phi_y\aeta-\eta_x\phi_x\aeta=B-\eta_x V$. 
Consequently, one has the identity
\begin{align*}
(G(\eta)\psi)^2&=B^2-2 \eta_x BV +\eta_x^2 V^2 \\
&=B^2-V^2-2\eta_x BV+(1+\eta_x^2)V^2.
\end{align*}
Now, in terms of $N(\eta)\psi$ (see \e{t90}), this gives
$$
(1+\eta_x^2)V^2=(G(\eta)\psi)^2+2 N(\eta)\psi.
$$
We next use that (cf \e{C6-bis})
$$
\int_{-L}^L \chi N(\eta)\psi\dx
=-\iint_{\Omega} \chi_x\phi_x\phi_y \dydx +\mez \int_{-L}^L \chi \phi_x^2\arrowvert_{y=-h}\dx,
$$
to infer that
\be\label{d10}
\begin{aligned}
\intL \chi (1+\eta_x^2)V^2 \dx &=
\intL \chi (G(\eta)\psi)^2\dx + \intL \chi\phi_x^2\arrowvert_{y=-h}\dx\\
&\quad-2\intL\int_{-h}^{\eta(t,x)}\chi_x\phi_x\phi_y\dydx.
\end{aligned}
\ee
To conclude the proof of \e{d11}, it remains to estimate 
$\intL \chi (1+\eta_x^2)\psi_x^2 \dx$ by means of $\intL \chi (1+\eta_x^2)V^2 \dx$. To do so, we start from the following identities 
(which easily follow from the definitions of $V,B$)
$$
V=\psi_x-\eta_x B, \quad B=G(\eta)\psi+\eta_x V,
$$
to obtain, using the assumption $\la \eta_x\ra\le 1$,
\begin{align*}
\psi_x^2&=(V+\eta_x B)^2\le 2V^2+2\eta_x^2 B^2\\
&\le 2V^2+4\eta_x^2(G(\eta)\psi)^2+4\eta_x^4V^2\\
&\le 4(1+\eta_x^2)V^2+4(G(\eta)\psi)^2.
\end{align*}
Then the wanted inequality \e{d11} follows from \e{d10}.
\end{proof}

Now, by combining the previous estimates, we conclude that
\begin{align*}
O+W+B-I&\le \left(C_1+K_1+K_2+\frac{1}{\eps_2\lambda}\sup_{[-L,L]} \chi+\frac{4\eps_3 L}{\lambda}+\frac{L}{2\eps_3\lambda}\right)\Hr(0)\\
&\quad +\left(2\eps_3 L-\frac{h}{4}\right)\intL \phi_x^2\arrowvert_{y=-h}\dxdt \\
&\quad + \left(\frac{\eps_1K_1}{2}+\frac{9\eps_2}{16}K_2 +4\eps_3 L\sup \la \chi_x\ra\right)\int \Hr(t)\dt,
\end{align*}
where
\begin{align*}
C_1&=\frac{\lambda}{2\eps_1}\sup_{[-L,L]} \chi(x) +\frac{3\lambda}{2g}\sup_{[-L,L]}(1-m_x(x))^2,\\
K_1&=\frac{6}{\kappa}\sup_{[-L,L]}m(x)^2 +\frac{4}{g}\sup_{[-L,L]}\left( \frac{5}{4}-\frac{m_x(x)}{2}\right)^2,
\end{align*}
and where $K_2$ is given by \e{Bures2}. 

Now, we fix $\lambda=1$ (we do not consider the problem of optimizing $\lambda$) and we fix 
$\eps_1,\eps_2$ and $\eps_3$ small enough, so that
$$
2\eps_3 L-\frac{h}{4}\le 0,\quad 
\frac{\eps_1K_1}{2}+\frac{9\eps_2}{16}K_2 +4\eps_3 L\sup \la \chi_x\ra\le \frac{1}{8}. 
$$
Then
\be\label{C13bis}
O+W+B -I \le K \Hr(0)+\frac{1}{8} \int_0^T\Hr(t)\dt,
\ee
where $K$ is defined by
\be\label{d20}
K\defn C_1+K_1+K_2+\frac{1}{\eps_2\lambda}\sup_{[-L,L]} \chi+\frac{4\eps_3 L}{\lambda}+\frac{L}{2\eps_3\lambda}.
\ee
Notice that $K$ depends only on $\kappa,g,L,h$.
Then we obtain \e{C11} with $C=8K$. 
This concludes the proof of Theorem~\ref{T:2}.

\appendix

\section{}

For the sake of completeness, we recall from \cite{A-Stab-WW} the proof of the identity~\e{C6}. 

Here the time variable is seen as a parameter and we drop it. 
We consider a real number $s>5/2$ and two $2L$-periodic functions $\eta,\psi\in H^s(\xT)$ which are even in $x$. 
We denote by $\phi$ the harmonic function defined by 
\begin{equation}\label{m1}
\left\{
\begin{aligned}
&\Delta_{x,y}\phi=0\quad\text{in }\Omega=\{(x,y)\in \xT\times \xR \,;\, -h<y<\eta(x)\},\\
&\phi(x,\eta(x)) = \psi(x),\\
&\phi_y(x,-h) =0.
\end{aligned}
\right.
\end{equation}
It follows that $\nabla_{x,y}\phi$ belongs to $C^1(\overline{\Omega})$ 
and that $\phi_x=0$ when $x=-L$ or $x=L$ (see Proposition $2.2$ in \cite{Boundary}). As a result 
one can take the trace of $\nabla_{x,y}\phi$ and define
\be\label{d5}
N(\eta)\psi=\mathcal{N} \big\arrowvert_{y=\eta}\quad\text{with}\quad
\mathcal{N}=
\mez\phi_x^2-\mez \phi_y^2+\eta_x \phi_x \phi_y.
\ee
\begin{lemm}
For any $C^1$ function $\mu=\mu(x)$ which is $2L$-periodic, there holds
\be\label{C6-bis}
\int_{-L}^L \mu N(\eta)\psi\dx
=-\iint_{\Omega} \mu_x\phi_x\phi_y \dydx +\mez \int_{-L}^L \mu \phi_x^2\arrowvert_{y=-h}\dx.
\ee
\end{lemm}
\begin{proof}
The proof relies on the following identity
$$
\partial_y \big( \phi_y^2-\phi_x^2\big)+2\partial_x \big(\phi_x\phi_y\big)=
2\phi_y \Delta_{x,y}\phi,
$$
which implies that, since $\phi$ is harmonic and $\partial_y \mu=0$,
$$
\partial_y \big( \mu \phi_y^2-\mu \phi_x^2\big)+2\partial_x \big(\mu \phi_x\phi_y\big)=
2\mu_x\phi_x\phi_y.
$$
We deduce that the vector field $X\colon \Omega\rightarrow \xR^2$ defined by $X=(-\mu \phi_x\phi_y;\frac{\mu}{2}\phi_x^2-\frac{\mu}{2}\phi_y^2)$ 
satisfies
$\cn_{x,y} \big(X\big)=-\mu_x\phi_x\phi_y$. 
Since $\nabla_{x,y}\phi$ belongs to $C^1(\overline{\Omega})$ and since one has 
the boundary conditions
$$
\phi_y\arrowvert_{y=-h}=\phi_x\arrowvert_{x=-L}=\phi_x\arrowvert_{x=L}=0,
$$
an application of the divergence theorem gives that
\begin{align*}
-\iint_\Omega \mu_x\phi_x\phi_y \dydx&=\iint_\Omega \cn_{x,y}X\dydx\\
&=\int_{\partial\Omega}X\cdot n \dsigma= \int_{-L}^L \mu \mathcal{N} \big\arrowvert_{y=\eta} \dx
-\mez\int_{-L}^L \mu \phi_x^2\arrowvert_{y=-h}\dx.
\end{align*}
This completes the proof.
\end{proof}

\section{}\label{Proof:Bures2}

Let us prove \e{Bures2}. The time variable is seen as a parameter and we drop it. 
As above, we introduce the harmonic extension of $\widetilde\psi$ defined by 
\begin{equation*}
\left\{
\begin{aligned}
&\Delta_{x,y}\widetilde\phi=0\quad\text{in }\Omega=\{(x,y)\in \xT\times \xR \,;\, -h<y<\eta(x)\},\\
&\widetilde\phi(x,\eta(x)) = \widetilde\psi(x),\\
&\widetilde\phi_y(x,-h) =0.
\end{aligned}
\right.
\end{equation*}
Since $\widetilde\psi=\psi-\langle\psi\rangle$ where $\langle\psi\rangle=\frac{1}{2L}\intL \psi(x)\dx$ is a constant, we have 
$\widetilde\phi=\phi-\langle\psi\rangle$ where $\phi$ is the harmonic extension of $\psi$. 

It will be useful to consider a diffeomorphism $\Sigma$ from the 
flat strip $\xT\times [-1,0]$ to the domain $\Omega$, of the form
$$
\Sigma\colon (x,z)\mapsto (x,\sigma(x,z))\quad 
\text{with }\sigma(x,z)=(1+z)\eta(x)+hz.
$$
Since $\eta\ge -h/2$ by assumption, we easily verify that $\Sigma$ is a diffeomorphism 
from $\xT\times [-1,0]$ to the domain $\Omega$. Then we set
$$
\widetilde\varphi(x,z)=\widetilde\phi(x,\sigma(x,z)).
$$
Recall that $\la \eta_x(t,x)\ra\le 1$ by assumption (cf Theorem~\ref{T:2}). 
Then, directly from the change of variables formula for integrals, we obtain that there exists a constant 
$C$ depending only on $h$, such that
$$
\lA \nabla_{x,z}\widetilde\varphi\rA_{L^2(\xT\times [-1,0])}^2\le C\blA \nabla_{x,y}\widetilde\phi\brA_{L^2(\Omega)}^2=
C\lA \nabla_{x,y}\phi\rA_{L^2(\Omega)}^2\le C \Hr.
$$
Consequently, to prove \e{Bures2} it is sufficient to prove that 
\be\label{m2}
\blA \widetilde\psi\brA_{L^2}^2\le C'\blA \nabla_{x,z}\widetilde\varphi\brA_{L^2(\xT\times [-1,0])}^2,
\ee
for some constant $C'$ depending only on $h,L$. To do so, we proceed in two steps. Firstly, 
we set
$$
c_n\defn \frac{1}{2L}\int_{-L}^L \widetilde{\psi}(x)e^{-\frac{\pi}{L}inx}\dx
$$
and use the fact that $\widetilde\psi$ is $2L$-periodic to obtain 
$$
\blA \widetilde\psi\brA_{L^2(\xT)}^2=2L\sum_{n\in \xZ}\la c_n\ra^2.
$$
By assumption, the mean of $\widetilde\psi$ vanishes and hence $c_0=0$. 
This immediately implies that
$$
\blA \widetilde\psi\brA_{L^2(\xT)}^2\le 2L\sum_{n\in \xZ} \la n\ra\la c_n\ra^2\le 
A \blA \la D_x\ra^{1/2}\widetilde\psi\brA_{L^2}^2,
$$
where the constant $A$ depends only on $L$ and 
$\la D_x\ra^{1/2}$ is the Fourier multiplier with symbol $\la\xi\ra^{1/2}$. 
Now, to bound $\blA \la D_x\ra^{1/2}\widetilde\psi\brA_{L^2}^2$ by 
$\blA \nabla_{x,z}\widetilde\varphi\brA_{L^2(\xT\times [-1,0])}^2$, we may 
proceed as in the proof 
of Proposition~$3.12$ in \cite{LannesLivre}.

\bibliographystyle{plain}

\vspace{3cm}

\begin{flushleft}
\noindent Thomas Alazard\\[1ex]
CNRS et Centre de Math\'ematiques et de Leurs Applications UMR 8536\\
\'Ecole Normale Sup\'erieure de Paris-Saclay\\
61 avenue du Pr\'esident Wilson\\
94235 Cachan Cedex
\end{flushleft}

\end{document}